\documentclass[12pt]{article}
\usepackage{epsfig,psfrag,amsmath,amssymb,latexsym}
\usepackage{amscd }
\usepackage{amsfonts}
\usepackage{graphicx}
\newtheorem{theorem}{Theorem}[section]

\newtheorem{lemma}{Lemma}[section]

\newenvironment{proof}[1][Proof]{\textbf{#1.} }{\hfill\rule{0.5em}{0.5em}}
\newenvironment{eqproof}[1][Proof]{\noindent\textbf{#1.} }{}

\newcommand{\PP}{\mathbb{P}}
\newcommand{\E}{\mathbb{E}}
\def\cale{\mathcal E}
\def\var{\mathop{\rm Var}}
\numberwithin{equation}{section}
\begin{document}
\title{On the Variance of the Optimal Alignments Score for
Binary Random Words
and an Asymmetric Scoring Function} 
\author{Christian Houdr\'e\thanks{School of Mathematics, Georgia
 Institute of Technology, Atlanta, GA 30332. Research supported
in part by the Simons Foundation Grant \#246283.
 ({\tt houdre@math.gatech.edu})}\; and Heinrich 
Matzinger\thanks{School of Mathematics, Georgia 
Institute of Technology, Atlanta, GA 30332.\hfill\break
({\tt matzi@math.gatech.edu})\hfill\break
Key Words: Optimal Alignments, Variance Bounds, Longest Common Subsequences,
Last Passage Percolation\hfill\break
 }}

\maketitle

\abstract We investigate the order of the variance of the 
optimal alignments score of two independent 
iid binary 
random words having the same length. The letters are equiprobable, but
the scoring function is such that one letter 
has a larger score than the other.  
In this setting, we prove that 
the order of variance is linear in the common length.
Optimal alignments constitute a generalization
of longest common subsequences, they can be represented
as optimal paths in a two-dimensional 
last passage percolation setting with 
dependent weights. 

\section{Introduction}
The problem under investigation in this article 
is similar to the corresponding problem for the longest common 
subsequence, which is, in fact, a special case of optimal 
alignment.  

To start, let us give a few definitions and present an example.
Consider two (non-random) strings  $x=x_1x_2\cdots x_n$ and $y=y_1y_2\cdots y_n$
 both  having length $n$ and written with an alphabet 
$\mathcal{A}$, and consider {\it alignments with gaps}
of  these two   strings. An alignment $\pi$ and, say,
if $n=5$ could, for example, be:
$$
\begin{array}{c|c|c|c|c|c|c}
x_1&   &x_2&x_3&x_4&   &x_5\\\hline
y_1&y_2&y_3&   &y_4&y_5& 
\end{array}
$$
Next,  consider a scoring function 
$s:\mathcal{A}\times\mathcal{A}\rightarrow \mathbb{R}$ which,
 in applications, measures the similarity/dissimilarity between
letters. The total score of an alignment is then
the sum of the individual scores of the aligned letter pairs, minus
a penalty proportional to the number of gaps. In the present example,
$S_\pi(x_1\cdots x_n;y_1\cdots y_n)$, the alignment score of $\pi$,
is given by
$$S_\pi(x_1\cdots x_5,y_1\cdots y_5)=s(x_1,y_1)+s(x_2,y_3)+s(x_4,y_4)-4q,$$
where $q>0$ is a {\it gap penalty}. 
An alignment which maximizes the alignment score of the strings
$x$ and $y$ is called an {\it optimal 
alignment} (OA). 
The score of an optimal alignment is an {\it optimal alignment
score} and is denoted by $S(x_1\cdots x_n;y_1\cdots y_n)$. 
The optimality of an  alignment
depends, of course, on the scoring function under consideration.
The length of the longest common subsequences (LCSs) of $x$ and $y$ can be viewed
as  an optimal alignment
score of $X$ and $Y$, with
$s(i,j)=\delta_{ij}$ and a $q=0$ gap-penalty. 

Throughout we consider  two binary iid random strings 
$X=X_1X_2\cdots X_n$ and $Y=Y_1Y_2\cdots Y_n$ independent of each
other. The optimal alignment score of $X$ and $Y$ is denoted by $L_n$
so that
$$L_n:=S(X_1\cdots X_n; Y_1\cdots Y_n),$$
and the main result of the present paper is that under a sufficiently 
large bias
in the scoring function $s$, then 
\begin{equation}
\label{VARLn}
\var L_n=\Theta(n),
\end{equation}
i.e., $\var L_n/n$ is upper and lower bounded by positive
constants independent of $n$.
For this
we take both the gap penalty $q$, $s(1,0)$ and $s(0,1)$ to all be 
equal to $0$, 
also assume that
$$\PP(X_i=1)=\PP(X_i=0)=\PP(Y_i=0)=\PP(Y_i=1)=\frac12,$$
and  further request that $s(1,1)$ be a given amount bigger than $s(0,0)$.
The order \eqref{VARLn} which was conjectured by Waterman 
\cite{Waterman-estimation},  had previously only been established
for some very special cases 
\cite{bonettolcs}, 
\cite{increasinglcs},
\cite{periodiclcs},  
like for one sequence non-random but periodic or having an
extra symbol or under an extra increasing requirement.  
The order \eqref{VARLn} has also been established
for a low-entropy case in \cite{LM} (see also \cite{HM}).
There, $L_n=LC_n$, the length of the 
longest common subsequences of $X$ and $Y$, 
and $\PP(X_i=1)=\PP(Y_i=1)$ is very close
to $0$. 
To the best of our knowledge, 
the present is the first time that for a non-low-entropy case,
this order is established for the optimal alignments score
of random strings. 

For the LCS-case, 
Steele \cite{Steele86}, proved, in particular, that 
${\rm Var }\; LC_n\leq n$, while 
Chv\'atal and Sankoff conjectured that 
${\rm Var }\; LC_n$ was of order $n^{2/3}$, which is the 
same order (when properly rescaled), as the one obtained 
by Baik, Deift and Johansson \cite{BaikDeiftJohansson99} 
in their much celebrated result on the longest increasing 
subsequence (LIS) of a uniform random permutation
of $\{1,2,\dots, n\}$.  
(See, Romik \cite{R} for an up-to-date account of this problem with
a complete bibliography.)
 Using a superadditivity argument, 
Chv\'atal and Sankoff \cite{Sankoff1} 
prove that 
$$\gamma^*:=\lim_{n\rightarrow\infty}\frac{\E\, LC_n}{n}$$
exists, where again $LC_n$ is the length of the LCS of 
two independent iid sequences of length $n$.  
To this day, the value of $\gamma^*>0$ 
is unknown even for binary sequences.  
 Arratia and Waterman \cite{watermanphase} derive 
large deviations for the fluctuations of $L_n$ on scales 
larger than $\sqrt{n}$.    Using first passage percolation
methods, Alexander \cite{Alexander} (see also \cite{ratelambda})
 proves that $\E\, L_n/n$ converges
to $\gamma^*$ at a rate of order at least $\sqrt{\log{n}/n}$.
The nature of the optimal alignment
has also been studied in \cite{multiplicity}.
Finally a CLT is obtained in \cite{HI} under a sublinear 
lower bound assumption
on $\var LC_n$ and this is further extended to the OA multi-sequences
framework in \cite{GHI}.

As already mentioned, our current 
OA-problem can be reformulated
as a LPP problem with dependent weights.
At this stage,
let us explain how this is so.
In LPP one considers an oriented graph $(\cale ,V)$
with a random coloring $w:\cale \rightarrow \mathbb{R}^+$ of the edges. 
An optimal path from $x$ to $y$ is then a path $(x_1,x_2,\ldots,x_m)\in V^m$
from $x$ to $y$ which maximizes the
total weight: $\sum_{i=1}^{m-1}w(x_ix_{i+1})$. (So, one requests
$x_1=x$, $x_m=y$ and $x_ix_{i+1}\in \cale$, for all $i=1,\ldots,m-1$.)
Now,
the optimal alignment score $L_n(S)$ is the weight of the heaviest path(s)
from $(0,0)$ to $(n,n)$. For this, take
the set of vertices $V$ to be $\mathbb{N}\times\mathbb{N}$
and the edges to always go one to the right or one up or diagonally
up to the next vertice. The weight for horizontal and vertical
edges is minus the gap penalty. For the edge $((i-1,j-1),(i,j))$
the weight is $s(X_i,Y_j)$. Within this setting, aligning a letter
with a gap corresponds to moving one unit vertically or horizontally.
Aligning $X_i$ with $Y_j$  corresponds to moving along the edge
$((i-1,j-1),(i,j))$. The optimal path then defines an optimal alignment:
for every edge $((i-1,j-1),(i,j))$ contained in the optimal path,
align $X_i$ with $Y_j$.

It is well known that the determination of 
the order of the variance for 
the LCS-problem is one of the main problems in computational biology, and to 
this day the generic problem remains open.  As  explained above, 
the optimal alignment problem can be reformulated 
 as an oriented last passage percolation (LPP)
problem with dependent weights.  For first passage percolation
(FPP) and LPP 
the exact order, in the general case, also remains unknown.  
FPP and LPP are part of a vast area of 
statistical physics \cite{spohn91} which is concerned with random growth models
for which physicists expect some universality properties.
More specifically one considers growth of a cluster where
material is being attached randomly on the surface of a nucleus. 
There are many fundamental
questions which have been open for decades, such as the universality of the 
fluctuation exponents.  
Physicists \cite{spohn91} have heuristic arguments
implying, among other results, that the fluctuations should behave like $n^{1/3}$ and the transversal
fluctuations should be of order $n^{2/3}$. 
(See the KPZ-conjecture in \cite{KPZ86}.)
 But, to this day, this has only 
been proven rigorously for some  special LPP models, like for the
Longest Increasing Subsequence 
of a uniform random permutation of $\{1,2,\dots, n\}$ \cite{BaikDeiftJohansson99} 
or LPP  on $\mathbb{Z}^2$ but with 
exponential or geometric waiting times 
\cite{Johansson2000}. For those cases, a Tracy-Widom limiting
distribution has been established. 
In the present
article we prove a totally different order. 

Optimal alignments have gained tremendous 
significance in both computational biology and computational
linguistics.  The reader will find in the standard references 
\cite{Backofen}, 
\cite{Durbin},  
\cite{Pevzner}, 
\cite{Sophie}, \cite{SK} and \cite{Water}
a general discussion of the relevance of 
string comparison, and related problems.

Optimal alignments and closely related methods are 
one of the main tools for identifying genes.  For example, with the help of
optimal alignments one searches for the location of gene X.  Assume 
that we have already identified this gene in one species, for example
the mouse.
Then probably a similar gene is present in the human genome. 
It is usually cheaper and more efficient to look
in the human genome for a substring similar to that mouse gene,
rather then to do again a full biological experimenting for human.
This implies that we have to be able to identify strings which are similar,
but where  some letters are missing. And this is exactly where 
optimal alignment scores are used: 
To evaluate the degree of similarity, 
the alignment score which is used is a log likelihood ratio
related to mutation probabilities.  A high alignment score indicates 
high similarity, but to understand the significance, a knowledge of
the variance of the optimal alignment score for 
the random model under consideration is needed.

Another  practical example of the use of  optimal alignments can be found in
 computational linguistics. A quite common task 
consists in determining, 
typically across languages, 
related pieces of texts or similar words.  
Often, one needs to identify pairs of translated words as this is an 
important step towards building electronic lexicons or a translation machine.  
Suppose we are given two texts: 
the first text is in English, while the second text is in German.  
The two texts are translations of one another and 
we try to build a computer program able to determine 
which words are translation of one another.  Many translation pairs 
exhibit great similarities.  Our program should be able to automatically
detect such similarities without prior knowledge of the languages.  
Let us present a concrete English/German example:    
Let $X$ be the English word $brother$ and let $Y$ be the German translation 
$Y=bruder$.  Looking at the two words, we immediately observe a great 
degree of similarity between them.  The computer should also 
be  able to detect this resemblance.  A first, unsophisticated 
method consists in writing one word beneath the other and to 
count the number of coinciding letters.  We then find:
$$\begin{array}{c|c|c|c|c|c|c}
b&r&o&t&h&e&r\\\hline
b&r&u&d&e&r&
\end{array}\,. $$
Two letters coincide: both words start with the letters $br$.  
This is not yet very conclusive.  For example the words
brag, bread, breast, bribe, bride, bring, broad, brute, all start 
with the letters $br$.   Hence this method of alignment 
is not very powerful to help discriminating unrelated pairs of words across 
similar languages.  A better solution consists 
in aligning the two words allowing for gaps while 
trying to obtain the maximum possible number of coinciding letters.  
With this method, the {\it optimal alignment} turns out to be:
$$\begin{array}{c|c|c|c|c|c|c|c|c}
b&r& & &o&t&h&e&r\\\hline
b&r&u&d& & & &e&r
\end{array}\,. $$
This time we get a sequence of four coinciding letters: 
$brer$.  Note that $brer$ is a common subsequence of $X$ and 
$Y$.  This means that the word $brer$ can be obtained from 
both $X$ as well as from $Y$ by only deleting letters.  
It turns out that $brer$ is the LCS of 
$X$ and $Y$.   
A further improvement consists in also allowing the alignment of 
similar letters.  For example, to give a score 
of one for identical letters, but a score of $1/2$ 
when the letters are only similar.  
Assuming  that $t$ and $d$ are similar, and that so are $o$ and $u$, 
we find the following optimal alignment:
\begin{equation}
\label{optalign}
\begin{array}{c|c|c|c|c|c|c}
b&r&o&t&h&e&r\\\hline
b&r&u&d& &e&r
\end{array}\,. 
\end{equation}
The score of the above alignment is 
$4+2(1/2)=5$, since four identical letters 
are aligned as well as two similar ones.  The 
score of the alignment \eqref{optalign} can be written as:
$$s(b,b)+s(r,r)+s(o,u)+s(t,d)+s(e,e)+s(r,r)=4+2(1/2)=5,$$
where $s(x,y)$ denotes the score obtained by aligning
the letter $x$ with the letter $y$.  Sometimes a gap penalty is also in use.
Note that the length of the LCS is equal to the optimal alignment 
score when the substitution matrix $s$ is taken to be the identity matrix 
and that a zero gap penalty is in force.

The determination of the order of the 
 variance for the optimal alignment score
of two random strings of length $n$  is what will be 
of concern to us in the rest of the text.

\section{The Main Result}

Throughout this paper $(X_i)_{i\ge 1}$ and $(Y_i)_{i\ge 1}$ 
are two independent iid sequences of Bernoulli random variables 
with parameter $1/2$.  
We also assume in everything that 
follows that the substitution matrix is 
such that: 
\begin{equation}
\label{score}
s(1,0)=s(0,1)=0\;\;{\rm and}\;\; s(0,0)=1, 
\end{equation}
while the gap penalty is taken equal to zero:
\begin{equation}
\label{q}
q=0.
\end{equation}
Consider the two strings of equal length $X:=X_1X_2\cdots X_n$ and
$Y:=Y_1Y_2\cdots Y_n$.  An {\it alignment} is a pair of increasing sequences
$(\pi,\nu)$ such that
$$1\leq\pi(1)<\pi(2)<\cdots<\pi(k)\leq n$$
and 
$$1\leq\nu(1)<\nu(2)<\cdots<\nu(k)\leq n,$$
where $\pi=\pi(1)\pi(2)\cdots\pi(k)$,
$\nu=\nu(1)\nu(2)\cdots\nu(k)$ and $k\leq n$.  
The score of the alignment $(\pi,\nu)$
is defined as
$$S_{(\pi,\nu)}:=\sum_{i=1}^{k} s(X(\pi(i)),Y(\nu(i))).$$
The {\it optimal score} $L_n$ is then defined as
$$L_n:=\max S_{(\pi,\nu)},$$
where the maximum is taken over all possible alignments
$(\pi,\nu)$.

The main result of the present paper asserts that if the score $s(1,1)$ is
large enough, then the variance of the optimal score is
of order $n$.  More precisely, 
\begin{theorem}
\label{mainresult}
Let $(X_i)_{i\ge 1}$ and $(Y_i)_{i\ge 1}$ be 
two independent sequences of iid Bernoulli random variables 
with parameter $1/2$.  Let the substitution matrices 
be such that $s(1,0)=s(0,1)=0$, $s(0,0)=1$ with moreover 
no gap penalty, i.e., $q =0$.  Then, 
there exist 
$s^*>0$ 
and 
$C>0$, 
such that if $s(1,1)\geq s^*$, then
\begin{equation}
{\rm Var} L_n\geq C n
\end{equation}
for all $n\in\mathbb{N}\backslash\{0\}$.
\end{theorem}
The above theorem implies that ${\rm Var} L_n=\Theta(n)$, i.e., 
$\var L_n$ is both lower and upper bounded by constants
(independent of $n$) times $n$, using the 
upper bound obtained, via the tensorization property of the 
variance, in \cite{Steele86}.  

The main idea developed in proving the above theorem, is to show that
changing the length of a randomly 
chosen block in $X$ has a tendency to increase the
score.  Since the number of blocks of a certain length 
has variance of order $\Theta(n)$ this, in turn, 
implies that ${\rm Var}L_n$ is of order 
$\Theta(n)$.

Let us present the heuristics for this proof and
let us now start with a formal definition. 

Let $0\leq i< j\leq n$.   The interval $[i,j]$ is called 
a {\it block of zeros} (resp.\ a {\it block of ones}) in $X$ if
$$X_{i}=X_{i+1}=\cdots=X_{j}=0, \quad (\mbox{resp.} = 1)$$
but $X_{i-1}=1$ (resp.\ 0) (or $i=1$) and $X_{j+1}=1$ (resp.\ 0) (or $j=n$).
The integer $j-i+1$ is then called the 
{\it length of the block} $[i,j]$.

For example, the string
$000011100000$ is made of three blocks, the first block consists
of four zeros and has length four,  the second block consists of 
three ones and has length three while the third block is a block of zeros 
of length five.

Next, let us described our transformation:  
Pick a block of zeros of length five at 
random among all 
the blocks of zeros of length five in $X$.  For this use the 
equiprobable distribution and this selection process 
is independent of $Y$.
Then, remove one zero from the chosen block.
The block becomes a block of length four. The next step is to add 
the zero just removed to a randomly chosen block of zeros of length one. 
This block then becomes a block of length two.  Again, 
to choose the block of length
one, use the equiprobable distribution on all blocks of zeros 
of length one in $X$.  
The string $X$ gets transformed in this
way into a new string having the same length.  This new string  is denoted
by $\tilde{X}$.

Once more, here is an example. Let
$X=0101000001000001$.  The string $X$
has two blocks of zeros of length one as well as two blocks
of zeros of length five.  Assuming that the second block 
of length five gets chosen and then so is the first block of length
one, we would obtain:
$\tilde{X}=0010100000100001$.

The optimal alignment score of $\tilde{X}$
with $Y$ is denoted by $\tilde{L}_n$. Hence:
$$\tilde{L}_n:=\max \sum_{i=1}^{k} s(\tilde{X}(\pi(i)),Y(\nu(i))),$$
where the maximum is taken over all the alignments $(\pi,\nu)$.
We show that, when
$s(1,1)$ is taken large enough, then  $\tilde{L}_n$ tends to be larger
than $L_n$. This is the content of the next theorem:

\begin{theorem}\label{maintheq}
Let $X=X_1\cdots X_n$ and $Y=Y_1\cdots Y_n$ and let both
\eqref{score} and \eqref{q} be satisfied. Let $\epsilon_1>0$
 and let $A^n$ be the event that the
following two conditions are satisfied:
\begin{equation}
\label{probcondA1}
\PP(\tilde{L}_n-L_n= 1\mid X,Y)\geq \frac{31}{128}
- \epsilon_1,
\end{equation}
\begin{equation}
\label{probcondA3}
\PP(\tilde{L}_n-L_n= -1\mid X,Y)\leq\frac{1}{32}
+\epsilon_1.  
\end{equation}
Then, there exist 
$s^*>0$
and 
$c_1>0$, 
both independent of $n$, 
such that if $s(1,1)\geq s^*$, then
\begin{equation}
\label{PAn}\PP(A^n)\geq 1-e^{-c_1n},
\end{equation}
for all $n\in\mathbb{N}\backslash\{0\}$.
\end{theorem}
The main idea behind the proof of this theorem is that
when $s(1,1)$ is large most of the ones will get matched with ones.
Let us explain the reason:
Think first of aligning as many $1$'s as possible
and nothing else.
This would lead to a score of about $s(1,1)n/2$, since
the proportion of $1$'s in both $X$ and $Y$ is about
$1/2$. Let $\delta>0$. Imagine an alignment leaving
 out a number $\delta n$
of $1$'s. (That number of left out $1$'s is to be compared to 
the maximum possible number of aligned $1$'s.)
That would represent ``a loss'' of $s(1,1)\delta n$
compared to aligning as many ones as possible. Say an alignment
is optimal while leaving out that
many $1$'s. To compensate
for this loss of $1$'s, we would need to align 
 at least $s(1,1)\delta n$ symbol $0$'s.
There are at most $n$ such symbols in each string $X$ and $Y$,
because both strings have length $n$.
This implies that:
$$\delta\leq \frac{1}{s(1,1)},$$
yielding a bound of about ${n/s(1,1)}$
on the maximum number of $1$'s, we could
possibly
leave out in an optimal alignment.
But the total number of $1$'s is about $n/2$.
Hence as a proportion of the total number of $1$'s,  
the proportion of left out $1$'s 
should not very much exceed  $2/s(1,1)$. Any number bigger than
that  by a fixed quantity,
will be an upper bound holding up to an exponentially small probability in $n$. 
This is to say that for any fixed number $q_0>2/s(1,1)$,  the event
that the optimal alignment 
aligns a proportion of $1$ at least equal to $q_0$ 
is exponentially small.

So, we are close to a situation where we match all the ones and match
as many zeros in between ones as possible. (Provided we chose $s(1,1)$ 
sufficiently large.) 
For  an alignment which tries to match all the ones, the blocks of zeros
between matched ones  are iid.  The distribution of 
the length of the blocks of zeros between matched ones
is approximately geometric with parameter $1/2$. 
(To simplify the exposition in this
section, take the strings $X$ and $Y$ to have random length so that
they both contain exactly $n/2$ ones. In this way,
when aligning all $1$'s we get iid parts.)

Let us look at another example of an alignment
obtained by matching all the ones.  Take
$X=1011000001$ and $Y=1001010100$.  When we match all the ones
and as many zeros as possible in between matched ones, we obtain
the following alignment:
$$
\begin{array}{cccccccccccccc}
1&0&\_ &1&\_&1&0&0 &0 &0&0&1 &\_&\_\\
1&0&0 &1&0 &1&0&\_&\_&\_&\_&1&0&0
\end{array}\,. 
$$ 
In our current example, the first block of zeros of $X$ has length $1$
and the second has length $5$.  Shortening the second block by removing a 
zero and adding it to the first block increases the alignment-score by one.
Indeed, since the second block of zeros of $X$ is matched with a  shorter
block, removing a zero does not reduce the score.  However now the first
block of zeros of $X$ is matched with a longer block of $Y$ so adding 
a zero increases the score by one unit.  When taking  a zero from the second
block and adding it to the first block of zeros of $X$, gives the following
new alignment:
$$
\begin{array}{ccccccccccccc}
1&0&0&1&\_&1&0&0 &0 &0&1&\_&\_\\
1&0&0 &1&0 &1&0&\_&\_&\_&1&0&0
\end{array}\,.
$$
This new alignment has score $4s(1,1)+3$,  which is one unit more than the score of 
the original alignment before modifications.

Let us consider yet another example. Take
the two strings of length $19$
$$x=1 0 1 01 000001 000001 0 1$$
and
$$y=1 001 01 000  1 00000 1001$$
and consider the alignmnent $\vec{\pi}$
which aligns all $1$'s and as many $0$'s as possible in between.
$$
\begin{array}{c|c|c|c|c| c|c|c|c|c| c|c|c|c|c| c|c|c|c|c|c| c|c|c}
x&&1 &0& &1 &0&1 &0&0&0&0&0&1 &0&0&0&0&0 &1 &0& &1\\\hline
y&&1 &0&0&1 &0&1 &0&0&0& & &1 &0&0&0&0&0&1 &0&0&1
\end{array}
$$
The total alignment score is then
given by 
$$S_{\vec{\pi}}(x,y)=6s(1,1)+1+1+3+5+1=6s(1,1)+11.$$
How large does $s(1,1)$ need to be for the above alignment which aligns
all $1$'s to be optimal? In the current case, if $s(1,1)>19$, then
we know that the optimal alignment cannot leave out any $1$'s, since
what one gets for leaving out just one pair 
of aligned $1$'s is losing more than $19$ points. The $0$'s cannot make
up for that since their maximal contribution would be $19$ points
if both strings would consist entirely of  $0$'s only.
(We assume exactly the same amout of $1$'s present in both strings
at this stage for simplification of the discussion.)

Now, let us apply our random transformation to $x$
to obtain the new string $\tilde{x}$. That is we chose 
a block of $0$'s of length $5$ at random in $x$ and remove a $0$,
which is then added to a randomly chosen block of 0's of length $1$.
In the present case, the expected increase in score for our alignment
$\vec{\pi}$ is:
\begin{equation}
\label{ES}
\E(S_{\vec{\pi}}(X,Y)-S_{\vec{\pi}}(\tilde{X},Y)|X=x,Y=y)=
\E(S_{\vec{\pi}}(\tilde{x},y)-S_{\vec{\pi}}(x,y))=\frac12(-1)+\frac{2}{3}>0
\end{equation}
Indeed, one of the two blocks of $0$'s of length $5$ is such that by
cutting one bit, the score gets reduced by one point,
while for the other block of $0$'s of length $5$,
the score does not change. This contributes
then  $-1/2$ to the expected change. Similarly two out of the three blocks
of $0$'s
of length $1$ in $x$ are such that adding a $0$ increase the score by one.
This gives a contribution of $2/3$ to the expectation.
Adding those two contributions then leads  to \eqref{ES}.\\

In general  the alignment which aligns
all ones is not an optimal alignment.  It is rather difficult to understand
how the optimal alignment looks macroscopically for long texts.  There are
extremely complicated
dependencies between all the parts of the optimal alignment. However, despite
the horrendously complicated macroscopic behavior,
when most ones are matched, the local distribution of the optimal alignment
is close to the alignment where we match all the ones. Moreover,
the alignment where we match all the ones consists of a sequence
of iid blocks of zeros between matched ones, provided we use a little
trick. That trick consists of taking the sequences $X$ and $Y$ having
random lengths. For this let $T_i$ be the 
position of the $i$-th $1$ in  $X_1,X_2,\ldots$:
$$X_1+X_2+\ldots+X_{T_i}=i$$
and
$$X_{T_i}=1.$$ 
Similarly, let $R_i$ be the position of the $i$-th $1$ in
the string $Y_1,Y_2,\ldots$. Now, consider the strings
with random length
$$X_1X_2\ldots X_{T_{n/2}},$$
and
$$Y_1Y_2\ldots Y_{R_{n/2}}.$$
Both these strings, have length about $n$ plus or minus a fluctuation of 
$\sqrt{n}$. They have the same amount of $1$'s.  Hence, when matching 
all the ones between these two strings, then  the blocks
of $0$'s in between  are exactly independent geometric random variables.
So, the proofs are done for such random 
length strings. One can then also alter the lengths to be 
$T_{n/2+\epsilon_1 n}$
and $R_{n/2+\epsilon_2 n}$ where $\epsilon_1,\epsilon_2$ are two small 
constants and specify which $1$ are not to be matched in advance
of drawing the strings.  There will be only a linear number 
of such $\epsilon_1$
and $\epsilon_2$ to be considered, so that the actual strings
$X$ and $Y$ admit one such representation.  
But the bias, for the strings 
with random lengths, will hold with very high probability.  
(From the independence of the parts between aligned $1$'s and exponential  
bounds, it holds with probability one minus an exponentially small quantity.) 
So, it will typically hold for all such strings at the same time, 
hence also for the original
$X$ and $Y$. In the Subsection~\ref{heiniheino}, these random
strings are used without mentioning it all the time 
and, there, the reader should think of 
$X$ and $Y$ as being such strings with random length.
In that section we provide a less detailed approach, 
which should nevertheless prove to be very useful 
to get an overall idea of why and how things work.  

\

For the alignment with all ones matched, 
we can compute the probability that the score increases by 
one when taking a zero from a block of $X$ of length five 
and adding it to a block of length one.
The blocks are chosen uniformly at random. The corresponding blocks in $Y$ are iid
and have approximately a geometric distribution.  Hence the probability to have the alignment
score increased by one is the probability that the chosen block of length  
five is matched with a block of shorter length times the probability that
the chosen block of length one is matched with a longer block.  Hence the
probability of an increase in score is approximately:
$$\PP(Z<5)\PP(Z>1)=\frac{31}{128},$$
where $Z$ is a geometric random variable with parameter $1/2$.

Similarly, the probability that the score decreases by one unit
is the probability that the block of length five is matched to a block
of length five or longer, times the probability that the block
of length one is matched to a block of length one or having no zero.  
Hence the probability that the score decreases by one unit is approximately:
$$\PP(Z\geq 5)\PP(Z\leq 1)=\frac{3}{128}.$$

Theorem \ref{maintheq} asserts that the probability for the optimal 
score to increase/decrease
through our bit transfer procedure is close to the probability
for the alignment matching all the ones. (For this we assume
that most ones are matched due to the score $s(1,1)$ being large.)  
Theorem \ref{maintheq} is proved in Section
4 and Section 5 using exponential estimates and some combinatorics.
In the next section, it is shown that Theorem \ref{maintheq} 
implies Theorem \ref{mainresult}. However, to begin, in the next subsection,
we compute a value for $s(1,1)$ which guarantees
that $\var L_n=\Theta(n)$.  There, for simplification, a slightly different 
notation than in the subsequent proofs is used.
So, the next subsection, can also be viewed as a quick overview
of the reasons for the method to work.

\subsection{How large does $s(1,1)$ need to be? Numeric Lower Bound.}
\label{heiniheino}
In this subsection, we determine a lower bound for $s(1,1)$ 
 to guarantee a positive bias in our random change, (which in terms 
implies the desired fluctuation order $\sqrt{n}$ times constant). More exactly,
 for all $n$, we want
\begin{equation}
\label{biasbias}
\PP(\E(\tilde{L}_n-L_n|X,Y)\geq \epsilon)\geq 1-\exp(-cn),
\end{equation}
to hold for some $\epsilon,c>0$ not
depending on $n$. The question is if an unrealistic large $s(1,1)$
is needed or not. For example in \cite{LM},
the authors consider binary i.i.d. sequences. To prove that
the variance of the LCS is $n$ times constant, they  assume
the  probability of one of the symbols
to be less than one over several billions!\\
In the current subsection, we show that $s(1,1)$ does not need to be 
in the billions!  Rather, a value in the thousands
is enough. And hence, the phenomena described in the current paper,
does not just appear in very extreme situations. Also, it should be
noted that our bound is not optimized. We are confident that with
more work one could get  smaller value. The approach
used here to determine  a numerical lower bound
for $s(1,1)$ is similar to the one used subsequently for the formal
proof. However, the exact method is different: it simply
turned out that what is most useful for a  formal proof,
is not the same way which allows to easily 
obtain a good lower bound. So, in this subsection we present
things in a less formal way. But, this subsection
has also the merit of allowing the reader, to understand
quickly what the main ideas of our proof are, without the burden
of all the technicalities in the later parts of this article.\\

To start with, consider an alignment which
aligns all $1$'s, instead of an optimal
alignment. For such an alignment the length of the blocks
of $0$'s in between $1$'s are approximately iid geometric variables with parameter
$1/2$. The expected increase due to a random change is then
typically about
\begin{equation}
\label{derbias}
\E(S_{\pi}(\tilde{X},Y)-S_{\pi}(X,Y)|X,Y)\approx \frac{1}{4}-\frac{1}{32}
=\frac7{32}>0,
\end{equation}
where $\pi$ denotes an alignment which aligns all $1$'s.

The advantage of the alignment which leaves out no $1$'s is
that the parts become independent.
Now, we would like to prove
the high likeliness of an inequality like \eqref{derbias}, but 
for the alignment $\pi$ being optimal instead of being the alignment
which leaves out no $1$'s.  
The problem with the 
optimal alignment is that it has complicated correlations between the different 
parts. This makes explicit calculations almost impossible. To circumvent
this difficulty, we will prove that a bias like in \eqref{derbias} holds
for all alignment $\pi$ which leave out only a small amount 
of
$1$'s. (The left out $1$'s referred to here, are the ones which get aligned with 
 gaps.) Then, we show that with high probability the optimal alignment 
is part of that collection of alignments. The advantage of
this approach is that we can work with alignments for which 
the parts inbetween aligned $1$'s are independent. For this
we determine before drawing $X$ and $Y$ which $1$'s are to be
aligned with gaps. Then, we align all other $1$'s with each other
and as many $0$'s as possible inbetween aligned $1$'s.  

Let us give an example. 
So, take for example $n=19$, and consider the strings $X$ and $Y$
to have random length and be taken so as to contain exactly six $1$'s each.
Next define an alignment of $X$ and $Y$ by specifying which
$1$'s are to be aligned with gaps. 
 This defines
 alignments  which boasts 
independent parts inbetween aligned $1$'s, allowing to use exponential inequalities.

Next, draw the strings  $X$ and $Y$ using an unbiased coin,
getting the strings $x$ and $y$:
$$x=1 0 1 01 000001 000001 0 1$$
and
$$y=1 001 01 000  1 00000 1001$$
So, align first those $1$'s aligned with gaps.
Then, all the remaining $1$'s get aligned with $1$'s. 
Finally inbetween aligned $1$'s align as many $0$'s
as possible. For example, request that the second $1$ of $x$ as well
as the fifth $1$ of $y$ get each aligned each with a gap.
This ``defines" then the alignment $\nu$ given by:
{\footnotesize
$$
\begin{array}{c|c|c|c|c| c|c|c|c|c| c|c|c|c|c| c|c|c|c|c|c| c|c|c|c|c| c|c|c}
x&&1 &0&1&0&1 &0&0&0&0&0&1 &0&0&0&0&0&1 &0& & & & & & & &1\\\hline
y&&1 &0& &0&1 &0& & & & &1 &0&0&0& & &1 &0&0&0&0&0&1&0&0&1
\end{array}
$$
}

\noindent
The alignment $\nu$ above can  also  be interpreted as follows: First,
leave out the $1$ to be aligned with gaps, that is concatenate
the blocks of $0$'s adjacent to those ``left out $1$'s".
Then align all remaining $1$'s with each other. 
And finally align a maximal number
of $0$'s from the new blocks of $0$'s in between aligned $1$'s. 
Some of the new blocks
will be original blocks and some have been obtained by concatenation.

In the example, the third  block of $0$'s of $x$
is aligned with the second block of $0$'s of $y$. Similarly,
the fourth block of $0$'s of $x$ is aligned with the third block
of $0$'s of $y$. So, the third and fourth blocks of $0$'s of $x$
are aligned ``one block onto one block'' from the original blocks.  
These blocks are {\it unproblematic} and are distributed like geometric variables with expectation $2$.

On the other hand, 
the first block of $0$'s of  $x$ in our example is different. It 
is ``concatenated''  with the second block of $0$'s of $x$. (Because
the $1$ inbetween these two blocks is aligned with a gap.)
Then this concatenated new block is aligned with the 
first block of $0$'s of $y$. This yields two aligned $0$'s.

This concatenation tends to decrease the likelihood that a random change
operated on a block of length one leads to an increase in score.
The reason being that
the added random bit may no longer increase the score:
 there might
already be enough $0$'s around due to the concatenation.

We call such a block of $0$'s of length $1$ which gets concatenated
(i.e., that is adjacent to a $1$ which gets aligned to a gap)
a {\it problematic block}. Which blocks are problematic of course
depends on the alignment under consideration.
Similarly a block of $0$'s of length $5$ is called {\it problematic}
if it gets aligned to a concatenated block of $y$. Finally,
in our example the last block of $0$'s of $x$ ``gets aligned with
a concatenated block of $y$''. This means that it gets aligned with 
the concatenation of the fourth and fifth block of $0$'s of $y$.
Now, this is not a problem: the probability for the random change
to lead to an increase is even bigger through that. Imagine
that we add a bit to such a block of $0$'s of length $1$ of $x$.
If that block is aligned with a concatenated block of $y$,
then the concatenated block tends to be longer than a single
block. Hence, there is a larger  probability that adding the bit increases 
the score. So, the last block of $0$'s in our example is
considered an unproblematic block. 

Now let $B_\nu$ be the event that the blocks of length $5$ and $1$ chosen
at random for our random transformation  $X\mapsto \tilde{X}$
are both non-problematic according to the alignment $\nu$. 
 In the present example, 
$$\PP(B_\nu|X=x,Y=y)=\frac{1}{3},$$
since in the alignment 
$\nu$ one out of three blocks of $X$
of length $1$ is not problematic, and both blocks of length $5$ in $X$
are not problematic. 

Next, we analyze the effect of the random change depending
on whether there are problematic blocks chosen by our random 
alteration. We find
\begin{align} \label{suess}
&\E(S_\nu(\tilde{X},Y)-S_\nu(X,Y)|X,Y)\\
&=\E(S_\nu(\tilde{X},Y)-S_\nu(X,Y)|X,Y,B_\nu)\PP(B_\nu|X,Y)\nonumber\\
&\qquad \qquad +
\E(S_\nu(\tilde{X},Y)-S_\nu(X,Y)|X,Y,B_\nu^c)\PP(B_\nu^c|X,Y)\nonumber\\
&\qquad\ge
\E(S_\nu(\tilde{X},Y)-S_\nu(X,Y)|X,Y,B_\nu)\PP(B_\nu|X,Y)-(1-\PP(B_\nu|X,Y)),
\nonumber
\end{align}
where the last inequality uses the fact that
the random change can decrease the score by at most one unit.
Now, with the non-problematic blocks, the expected increase
of the random change can be treated
in the same way as for the alignment which leaves out no $1$'s.
That is the expected increase is about $7/32$.
But, again, we want to prove a bias not for a single alignment but
for all the alignments which leave out less than a proportion $q_0$ of $1$'s.
There is an exponential number of such alignments. Consider the maximum
$$\max_\mu \E(S_\mu(\tilde{X},Y)-S_\mu(X,Y)|X,Y,B_\mu),
$$
where $\mu$ ranges over all alignments of $X$ and $Y$ which leave out less
than a proportion $q_0$ of $1$'s. For one alignment the expected change
should typically be $7/32$. Since, there are exponentially
many such alignments for which we wish to bound simultaneously
the expected increase, this number has to be lowered.
So, instead we take $1/5-4/32={3/40}$.
So, we consider the event $EVENT1$ that for all
alignments $\mu$ of $X$ and $Y$ which leave out no more than a proportion
$q_0$ of $1$'s, we have
$$\E(S_\mu(\tilde{X},Y)-S_\mu(X,Y)|X,Y,B_\mu)\geq \frac{1}{5}-\frac{4}{32}=
\frac{3}{40}.$$
Now, every $1$ aligned with a gap can create at most two problematic blocks
of $0$'s.  There are about $n/8$ blocks of $0$'s in $x$
of unit-length. Assume at first this number to be exact. Hence, if there
are no more than a proportion  $q_0$ of $1$'s aligned with gaps,
this yields at most a proportion of $8q_0$ problematic blocks of unit-length
$1$. (Counted among all blocks of length $1$ in $X$.) Similarly, we find
a proportion of problematic blocks of length $5$ not exceeding
$128q_0$. (Assuming the number of blocks of $0$'s of length
$5$ to be exactly $(n/2)(1/64)=(n/128)$.)
This then yields that if $\mu$ leaves out less than
a proportion $q_0$ of $1$, then:
\begin{equation}
\label{cafeone}
\PP(B_\mu|X,Y)\geq (1-8q_0)(1-128q_0)> 1-136q_0.
\end{equation}
This bound was obtained by assuming the number of blocks
of length $1$ and $5$ to be exactly equal to 
their expectation. In reality, we can only guarantee
that, with high probability, these numbers are close to it, up to a small
fixed error term. So, taking the bound on the rightmost side of \eqref{cafeone},
solves this problem. Hence,
\begin{equation}
\label{B}
\PP(B_\mu|X,Y)\geq 1-136q_0
\end{equation}
is therefore so is likely to hold for all alignments $\mu$ leaving out
no more than a proportion $q_0$ of $1$'s, and the likeliness of
inequality \eqref{B} to hold, up to an exponential
small quantity in $n$.

Assume next that the event $EVENT1$ holds. Using
\eqref{B} in \eqref{suess}, then
for any alignment $\mu$ leaving out less than a proportion
$q_0$ of $1$'s, it follows that
\begin{equation}
\label{q_0}
\E(S_\mu(\tilde{X},Y)-S_\mu(X,Y)|X,Y) \geq 
\frac{3}{40}(1-136 q_0)+ 136q_0
\end{equation}
In other words, in order to guarantee  a positive bias due to our random
change, 
$q_0$ needs to be such that:
\begin{equation}
\label{lll}
\frac{3}{40}(1-136q_0)-136q_0>0.
\end{equation}
This guarantees the bias we want, but only if the event
$EVENT1$ holds. Note that the inequality \eqref{lll}, holds
for 
\begin{equation}
\label{SoM}q_0\leq 0.0005299589
\end{equation} 
Now, for the event $EVENT1$ to hold, we want, 
for every alignment with no more than  a proportion
of $q_0$ gaps, the following two things to be true:
\begin{enumerate}
\item{}Among the non-problematic blocks of $0$'s of 
length $1$ of $X$, at least $1/5$ of them need
to be aligned with a block of length at least $2$.
Let the number of non-problematic blocks of length $1$ 
be denoted by $m$. For each such block, we have the event
that it is aligned with a block of length at least $2$.
These events are iid holding with probability $1/4$. Hence,
the total number of these non-problematic blocks, of length $1$
of $x$, satisfying this condition is a binomial variable
(when conditioning on $m$).
Hence, the probability not to have at least $1/5$ of these
blocks aligned with blocks of length bigger or equal to $2$ is
\begin{align}
\label{probabilitybound}
\PP\left(Binomial\left(m,\frac14\right) < \frac{1}{5} m\right)&\approx
\PP\left(Binomial\left(\frac n4,\frac14\right)\leq \frac{5n}{100}\right)\\
&\leq 
\left(\frac{\left(\frac14\right)^{\frac15}\left(\frac34\right)^{\frac45}}
{\left(\frac15\right)^{\frac15}\left(\frac45\right)^{\frac45}}\right)^{\frac n4}
\nonumber
\end{align}
where the last inequality was obtained using the fact that
$m\approx n/4$ and Lemma~\ref{B-lem}.
The probability bound given above in \eqref{probabilitybound} is for one
alignment. Now, we have exponentially many alignments
with less than a proportion $q_0$ of left out $1$'s. 
To be sure our bound still works for all those alignments
we need to make sure that the bound on the right side of \eqref{probabilitybound}
is much smaller than one over the number of such alignment.
For the number of alignments which leave out a proportion $q_0$ of $1$'s,
we use the bound in Lemma~\ref{A}. This then yields:
\begin{equation}
\label{SoM2}0.998251=
\left(\frac{\left(\frac14\right)^{\frac15}\left(\frac34\right)^{\frac45}}
{\left(\frac15\right)^{\frac15}\left(\frac45\right)^{\frac45}}\right)^{\frac14}
<q_0^{q_0}(1-q_0)^{1-q_0}
\end{equation}
This is the second condition on $q_0$ to ensure the bias.
\item{}Finally for the event $EVENT1$ to hold, the proportion of
 non-problematic blocks of $0$'s of length $5$ which are aligned with
a block of length greater or equal to $5$ needs to be less than  
${4/32}$. The probability for such an unproblematic block
to be aligned with a block of length at least $5$ is ${1/32}$,
so this total number is like a binomial random variable with parameter 
$m$ and $p={1/32}$. Here $m$ is the total number of blocks
of $0$ in $x$ of length $5$. This number is approximately
$m={n/64}$. Hence the probability of getting to many unproblematic
blocks of length $5$ aligned with blocks of length at least $5$ is about
\begin{equation}
\label{far-right}
\PP\left(Binomial\left(\frac{n}{64},\frac{1}{32}\right)\geq
 \frac{n}{64} \frac{4}{32}\right)
\leq \left(\frac{\left(\frac{1}{32}\right)^{\frac{4}{32}}
\left(\frac{31}{32}\right)^{\frac{28}{32}}}
{\left(\frac{4}{32}\right)^{\frac{4}{32}}
\left(\frac{28}{32}\right)^{\frac{28}{32}}}\right)^{\frac{n}{64}}
\end{equation}
Again this probability bound needs to be effective for all 
alignments which leave out less than a proportion of $q_0$ $1$'s.
So the probability bound on the right of \eqref{far-right}
needs to be strictly less than the number of such alignment.
With the help of Lemma~\ref{A} and Lemma~\ref{B-lem}, this leads to a third condition
on $q_0$:
\begin{equation}
\label{SoM3}0.9986848=\left(\frac{\left(\frac{1}{32}\right)^{\frac{4}{32}}
\left(\frac{31}{32}\right)^{\frac{28}{32}}}
{\left(\frac{4}{32}\right)^{\frac{4}{32}}
\left(\frac{28}{32}\right)^{\frac{28}{32}}}\right)^{\frac{1}{64}}
\leq q_0^{q_0}(1-q_0)^{1-q_0}
\end{equation}
\end{enumerate}

The question of {\bf how small $q_0$ needs to be
 to ensure that our results hold} can now be answered. For this,
all three equations \eqref{SoM}, \eqref{SoM2} and \eqref{SoM3} 
need to be simultaneously
satisfied. This is the case as soon as
$$q_0\leq 0.0001.$$
Since, we have seen that for an optimal alignment the proportion $q_0$ of 
leftout $1$'s is at most $1/s(1,1)$, the {\bf condition on s(1,1)
to guarantee} $\var L_n=\Theta(n)$ is then
$s(1,1)\geq 10^4$.

Let us quickly give the lemmas which
were used for the calculations above.
\begin{lemma}
\label{A}
Let $q\in[0,1]$. Consider the number of binary strings of length $n$
having $qn$ times $1$ appearing and the number of  $0$'s being
$(1-q)n$. Then, 
$$\binom{n}{qn}\leq \left(\left(\frac{1}{q}\right)^q
\left(\frac{1}{1-q}\right)^{1-q}\right)^n$$
\end{lemma} 
\begin{eqproof}Let $Z$ be a binomial variable with parameter $n$
and $q$. The probability that $Z$ equals $qn$ is less than one,
and hence:
$$\PP(Z=qn)=\binom{n}{qn} q^{qn}(1-q)^{(1-q)n}\leq 1$$
and therefore
$$\binom{n}{qn}\leq \left(\left(\frac{1}{q}\right)^q
\left(\frac{1}{1-q}\right)^{1-q}\right)^n.
\eqno{\rule{0.5em}{0.5em}}$$
\end{eqproof}

The next lemma gives an exponential bound for binomial variables.

\begin{lemma}
\label{B-lem}
Let $p,q\in[0,1]$ with $p\neq q$.
Let $Z$ be a binomial variable with parameter $p$ and $n$.
Then, if $q>p$
$$\PP(Z\geq nq)\leq 
\left(\frac{p^q(1-p)^{1-q}}{q^q(1-q)^{1-q}}\right)^n$$
while if $q<p$
$$\PP(Z\leq  nq)\leq 
\left(\frac{p^q(1-p)^{1-q}}{q^q(1-q)^{1-q}}\right)^n.$$
\end{lemma} 

\begin{eqproof}Assume $q>p$, the other case being similar.
We bound $\PP(Z=nq)$, but the same type of proof works
for $\PP(Z\geq nq)$.
So, $\PP_p$ (resp.\ $\PP_q$) refers to the 
probability measure when $Z$ has parameter
$p$ (resp.\ $q$) and $n$.
Then,
\begin{align*}
\PP_p(Z=nq)=
\binom{n}{nq} p^{nq}(1-p)^{n(1-q)}&=
\binom{n}{nq} q^{nq}(1-q)^{n(1-q)}
 \frac{p^{nq}(1-p)^{n(1-q)}}{q^{nq}(1-q)^{n(1-q)}}\\
&=
\PP_q(Z=nq)\left(\frac{p^q(1-p)^{1-q}}{q^q(1-q)^{1-q}}\right)^n
\end{align*}
since probabilities are always less or equal to one,
we find
$$\PP_p(Z\geq nq)\leq 
\left(\frac{p^q (1-p)^{1-q}}{q^q(1-q)^{1-q}}\right)^n.
\eqno{\rule{0.5em}{0.5em}}$$
\end{eqproof}

\section{Proof of the Main Theorem}
The purpose of this section is to prove that Theorem~\ref{maintheq}
implies our main Theorem~\ref{mainresult}.
To do so, we first need a few definitions.

Let $N_i$ be the number of blocks
of zeros of length $i$ in the string $X=X_1X_2\cdots X_n$, and 
let $\vec{N}$ be the vector $\vec{N}=(N_1,N_2,N_4,N_5)$.  
Let $H_a=\{(n_1,0,n_4,n_5):  n_1,n_4,n_5\in\mathbb{N}\}$
and let $H_b=\{(n_1,n_2,0,n_5): n_1,n_2,n_5\in\mathbb{N}\}$.
Let $\vec{e}:=(-1,1,1,-1)$.  It is easy to check that every 
integer vector $\vec{n}=(n_1,n_2,n_4,n_5)\in\mathbb{N}^4$ 
can be written in a unique way as a vector of $H_a\cup H_b$
plus a term $n\vec{e}$, where $n\in \mathbb{N}$.
For every $\vec{n}\in H_a\cup H_b$, define a sequence of
random strings
$$X(\vec{n}),X(\vec{n}+\vec{e}),X(\vec{n}+2\vec{e}),\ldots, X(\vec{n}+k\vec{e}),\ldots$$
as follows: The sequence is defined by induction on $k$:
\begin{itemize}
\item First, $X(\vec{n})$ is a binary random string of length $n$
having distribution $\mathcal{L}(X|\vec{N}=\vec{n})$. We also ask
that $X(\vec{n})$ be independent of $Y$.
\item Once, $X(\vec{n}+k\vec{e})$ is defined, let $X(\vec{n}+(k+1)\vec{e})$
be the string obtained by taking one zero from a block of
length five and adding it to a block of length $1$.
For this we draw uniformly a block among all the blocks
of zeros of $X(\vec{n}+k\vec{e})$ of length five and then
we pick again uniformly 
another block of zeros of length one in $X(\vec{n}+k\vec{e})$ 
and turn it into a block of length two.
\end{itemize}

Let us look, once again, at a numerical example. Take $n=16$ and
$\vec{n}=(2,0,0,2)$.  Hence the string $X(\vec{n})$ in this
case, is a binary random string of length $16$.  Its distribution
is the conditional distribution of an iid string, given
that there are:
\begin{itemize}
\item two blocks of zeros of length one,
\item no blocks of zeros of length two or four,
\item two blocks of zeros of length five.
\end{itemize}
We generate the string $X((2,0,0,2))$ with this distribution.
We could, for example, get
$$X((2,0,0,2))=0101000001100000.$$
To obtain the next string, that is
$X((1,1,1,1))$ we pick one of the blocks of length five at random
and turn it into a block of length four.  Then, we add a zero 
to one of the block of length one.  There are two blocks
of length five, so they have both equal probability 1/2 to get chosen.
Assume that the randomly chosen block of length five is the second one, and also
assume that the randomly chosen block of length one turns out to
be the first block.  Then, we would have that
$$X((1,1,1,1))=0010100000110000.$$

Let $L(\vec{n})$ be the optimal alignment score of
$X(\vec{n})$ and $Y$.
The next lemma asserts that for any $\vec{m}\in\mathbb{N}^4$, such that
$$\PP(\vec{N}=\vec{m})>0,$$
$L(\vec{m})$ has the same distribution of $L_n$ conditional
on $\vec{N}=\vec{m}$.
\begin{lemma}
Let $\vec{m}\in\mathbb{N}^4$ be such that
$$\PP(\vec{N}=\vec{m})>0.$$
Then, $L(\vec{m})$ has distribution
$$\mathcal{L}(L_n\mid \vec{N}=\vec{m}).$$
\end{lemma}
\begin{proof}We need to prove that
for any $\vec{m}\in\mathbb{N}^4$ with $\PP(\vec{N}=\vec{m})>0$,
$$\mathcal{L}(X(\vec{m}))=\mathcal{L}(X\mid \vec{N}=\vec{m}).$$
The proof is by induction on $k$.
By definition, 
$X(\vec{n})$ has distribution $\mathcal{L}(X|\vec{N}=\vec{n})$
when $\vec{n}\in H_a\cup H_b$. 
 The next step is to prove that when
$X(\vec{n}+k\vec{e})$ has distribution:
$$\mathcal{L}\left(X\mid \vec{N}=\vec{n}+k\vec{e}\right),$$
then   $X(\vec{n}+(k+1)\vec{e})$ must  have distribution
$$\mathcal{L}\left(X\mid \vec{N}=\vec{n}+(k+1)\vec{e}\right).$$
Let $\vec{m}\in \mathbb{N}^4$ be such that
$\PP(\vec{N}=\vec{m})\neq 0$. The distribution 
$\mathcal{L}(X|\vec{N}=\vec{m})$ can be characterized as
the uniform distribution on the set of strings $x$
for which
$$\vec{N}(x)=\vec{m}.$$
Here $\vec{N}(x)\in \mathbb{N}^4$ is the 
vector $(n_1^x,n_2^x,n_4^x,n_5^x)$,
where $n_i^x$ denotes the number of blocks of zeros
of length $i$ in the string $x$. 
Let 
$\vec{n}^*:=\vec{n}+k\vec{e}$, and then set
$$\vec{n}^*=(n_1^*,n_2^*,n_4^*,n_5^*).$$
So we only need to prove that all the possible realizations
of $X(\vec{n}^*+\vec{e})$ are equiprobable 
and that for every $x$ in the set 
\begin{equation}
\label{set}
\left\{x\in\{0,1\}^n:\vec{N}(x)=\vec{n}^*+\vec{e}\right\},
\end{equation}
 the probability 
$$\PP(X(\vec{n}^*+\vec{e})=x),$$
is non-zero.

Any string $x$,
such that $\vec{N}(x)=\vec{n}^*+\vec{e}$
has a non-zero number of blocks of zeros
of length two and  of
length four. Take any of these blocks of length four
and any block of length two. Let $y$ be the string
obtained by reducing the chosen block of length two by one bit
and adding this bit to the block of length four.
By the induction hypothesis, 
the probability for $X(\vec{n}^*)$
to be equal to $y$ is non-zero. The conditional probability
that $X(\vec{n}^*+\vec{e})=x$
given $X(\vec{n}^*)=y$, 
is equal to $1/(n_5^*n^*_1)$ which is different from zero.
Therefore,
$$\PP(X(\vec{n}^*+\vec{e})=x)\geq \frac{\PP(X(\vec{n}^*)=y)}{n_5^*n_1^*}>0.$$
It remains to prove that for all $x$ such that
$$\vec{N}(x)=\vec{n}^*+\vec{e},$$
$\PP(X(\vec{n}^*+\vec{e})=x)$ does not depend on $x$.
First,
\begin{equation}\label{y}
\PP(X(\vec{n}^*+\vec{e})=x)=\sum_y
\PP(X(\vec{n}^*+\vec{e})=x\mid X(\vec{n}^*)=y)
\PP(X(\vec{n}^*)=y),
\end{equation}
where the sum is taken over all strings $y$
which can be transformed into the string $x$
by taking a bit from a block of zeros of length five
and adding it to a block of zeros of length one.  For each pair consisting
of a block of zeros of length four and a block of zeros
of length two of $x$ there is such a $y$.  Hence the number of
terms in the sum on the right side of \eqref{y}, is equal to 
$$n_2^x n_4^x=(n^*_2-1)(n^*_4-1).$$
By the induction hypothesis, for a string $y$ such that
$\vec{N}(y)=\vec{n}^*$, 
$\PP(X(\vec{n}^*)=y)$ does not depend on $y$.
The probability for
$y$ to get transformed into $x$, is then equal to:
$$\PP(X(\vec{n}^*+\vec{e})=x\mid X(\vec{n}^*)=y)=\frac{1}{n^*_1n^*_5}.$$
Hence \eqref{y} becomes:
\begin{equation}\label{y2}
\PP(X(\vec{n}^*+\vec{e})=x)=
\PP(X(\vec{n}^*)=y)
\frac{(n_2^*+1)(n_4^*+1)}{n_1^*n_5^*},
\end{equation}
where $y$ is any string such that $\vec{N}(y)=\vec{n}^*$.
As already mentioned, by 
the induction hypothesis, for such a $y$, 
$\PP(X(\vec{n}^*)=y)$ does not depend on $y$. Hence,
the expression on the right side of  \eqref{y2} does not depend
on $x$. This finishes the proof.
 \end{proof}

We assume that the random variables $L(\vec{n})$
were constructed so that they are independent
of $\vec{N}$. Then
$L(\vec{N})$ has same distribution as $L_n$.
Hence,
\begin{equation}\label{Ln=LN}
{\rm Var } L_n={\rm Var } L(\vec{N}).
\end{equation}
We mentioned that any vector of $\mathbb{N}^4$
can be represented in a unique way as a sum:
$\vec{n}:=m_1+k\vec{e}$ where $m_1\in H_a\cup H_b$ and
$k\in\mathbb{N}$.
Let $\vec{M}\in H_a\cup H_b$ and $M_1\in \mathbb{N}$ be
defined by the equation:
\begin{equation}
\label{M}\vec{N}=\vec{M}+M_1\vec{e}.
\end{equation}
We also write
$L(M_1,\vec{M})$ for $L(\vec{N})$.
Let $H^n=
(H_a\cup H_b) \cap B(0,\sqrt{n})$, where $B(0,\sqrt{n})$
is the ball of radius $\sqrt{n}$ centered at $0$.
Let $H_e^n:=\{k\vec{e}:k\in [-\sqrt{n},\sqrt{n}]\}$.
Let $\vec{\mu}:=(\mu_1,\mu_2,\mu_4,\mu_5)$
where 
$$\mu_i:=\E N_i,\quad i=1,2,4,5 ,$$ 
and finally let 
\begin{equation}
\label{defIn}I^n:=C H^n+C H_e^n+\vec{\mu},
\end{equation}
where $C>0$ is a constant independent of $n$,
chosen large enough so that the following inclusion
\begin{equation}\label{IC}
\vec{\mu}+[-\sqrt{55n/4},\sqrt{55n/4}]^4
\subset I^n,
\end{equation}
is satisfied.
Let $E^n_{\rm slope}$ be the event that
for all $\vec{n}_a,\vec{n}_b\in I^n$ such that 
$\vec{n}_b-\vec{n}_a=k\vec{e}$, with $k\geq n^{1/10}$,
we have that
$$L(\vec{n}_b)-L(\vec{n}_a)\geq \frac1{100} |\vec{n}_b-\vec{n}_a|.$$
For any random variables $V$ and $W$, as usual
${\rm Var}(V\mid W) =\E(V^2\mid W)-(\E(V\mid W))^2$, hence
$${\rm Var}\, V=\E({\rm Var}(V\mid W))+{\rm Var} (\E(V\mid W)) ,$$
and so
$${\rm Var }\, V\geq \E { \rm Var }\,(V\mid W).$$
This gives in our case,
\begin{equation}
\label{beginning}
{\rm Var }\,L(\vec{N})\geq \E
({\rm Var}(L(M_1,\vec{M})\mid L(\cdot),\vec{M})).
\end{equation}
Since the variance is positive, we find
\begin{align}\label{alin1}
 \E\, &{\rm Var }(L(M_1,\vec{M})\mid L(\cdot ),\vec{M})\\
&\label{alin2}
\quad \geq \E({\rm Var }(L(M_1,\vec{M})|L(\cdot),\vec{M})\Big|
\vec{N}\in I^n, 
E^n_{\rm slope})
\PP(\vec{N}\in I^n) \PP(E^n_{\rm slope}).
\end{align}
Note that when we condition on $\vec{M}$ and hold $\vec{M}$ fixed,
then $L(M_1,\vec{m})$ becomes a function of the one dimensional
variable $M_1$. When the event $E^n_{\rm slope}$ holds,
that function has ``positive slope on the scale $n^{1/10}$".
We are therefore interested in the variance of a non-random function
of a random variable. For this, 
assume that $f:\mathbb{R}\rightarrow\mathbb{R}$ is such that
$f'(x)>c$, for all $x\in\mathbb{R}$ (where $c>0$ is a 
constant not depending on $x$). Then,
for any random variable $Y$, we have
\begin{equation}
\label{VarfY}
{\rm Var}\,f(Y)\geq c^2 {\rm Var}\,Y.
\end{equation}
Hence, if the map $x\mapsto L(x,\vec{M})$ would have partial derivative
along $x$ everywhere larger than
$c>0$, is would follow that ${\rm Var}\, L(M_1,\vec{m})$ is larger than 
$c^2 {\rm Var}( M_1\mid \vec{M})$.
Typically, the integer map $x\mapsto L(x,\vec{M})$ 
is not strictly
increasing in the direction of $\vec{e}$.
 But it is likely that the event $E^n_{\rm slope}$
holds, and hence there is a linear increase in the direction of
$\vec{e}$  every $n^{1/10}$ points.
We are next going to formulate a lemma which provides a modification
of \eqref{VarfY}, for when the map $f(\omega)$ does not increase
for every  $k$, but has a tendency to increase on some scale:

\begin{lemma}
\label{discretebound}
Let $c,m>0$ be two constants. Let $f:\mathbb{Z}\rightarrow 
\mathbb{Z}$ be a
non decreasing function such that:
\begin{itemize}
\item{} for all $i<j$:
\begin{equation}
\label{discretecondition1}
f(j)-f(i)\leq (j-i)
\end{equation}
\item{} for all $i,j$ such that $i+m \leq j$:
\begin{equation}
\label{discretecondition2}
f(j)-f(i)\geq c (j-i).
\end{equation}
\end{itemize}
Let $T$ be an integer-valued random variable such that $\E(|f(T)|)<+\infty$.
Then: 
\begin{equation}\label{fundi}
{\rm Var}\, f(T)\geq 
c^2\left(1-\frac{2m}{c\sqrt{{\rm Var}\, T}}
\right){\rm Var}\, T.
\end{equation}
\end{lemma}
We can now use Lemma~\ref{discretebound}
to give a lower bound on ${\rm Var}(L(M_1,\vec{M})\mid L(\cdot ),\vec{M})$
when $E^n_{\rm slope}$ holds and conditioning
on $\vec{N}\in I^n$. When $E^n_{\rm slope}$ holds,
the integer map 
$$x\mapsto L(x,\vec{M}),$$
restricted to $\{x: x\vec{e}+\vec{M}\in I^n\}$
satisfies the conditions of Lemma~\ref{discretebound}
with $m=n^{1/10}$ and $c=1/100$.  Hence, we obtain
that conditional on $\vec{N}\in I^n$, and when $E^n_{\rm slope}$
holds, then
\begin{align}\label{fundi2}
&{\rm Var}(L(M_1,\vec{M})\mid L(\cdot ),\vec{M})\\
&\quad \geq 
\frac1{10^4}\left(1-\frac{2n^{1/10}}{10^{-2}\sqrt{{\rm Var}(M_1\mid \vec{M},\vec{N}\in I^n)}}
\right){\rm Var}(M_1\mid \vec{M},\vec{N}\in I^n).\nonumber
\end{align}
Next we need the following lemma:
\begin{lemma}
\label{3.3}
Let $W$ be a random variable taking its values in $\mathbb{Z}$.
Let $J$ be an interval of diameter at least
$3\ln2\sqrt{n}/\kappa$, $n\ge 1$, for some constant $\kappa>0$,
independent of $n$, such that
$\PP(W\in J)=1$. Assume that
\begin{equation}
\label{Wkappa}
\PP(W=i+1)\geq \PP(W=i)\left(1-\frac{\kappa}{\sqrt{n}}\right),
\end{equation}
and 
\begin{equation}
\label{Wkappa1}
\PP(W=i)\geq \PP(W=i+1)\left(1-\frac{\kappa}{\sqrt{n}}\right),
\end{equation}
for all $i,i+1\in J$.
Then,
\begin{equation}
\label{Varkappa2}
{\rm Var}\, W\geq n\frac{(\ln 2)^2}{16\kappa^2}
\end{equation}
for every $n$ large enough.
\end{lemma}
\begin{eqproof}
Let $I$ be the interval
$$I:=\left[\E(W)-\sqrt{n}\ln2/(2\kappa),\E(W)+\sqrt{n}\ln2/(2\kappa)\right].$$
Let $I_r$ (resp.\ $I_\ell$) denote the interval of length $\ln2\sqrt{n}/\kappa$
directly adjacent to the right of $I$ (resp.\ to the left of $I$).
Then either $I_r\subset J$ or $I_\ell\subset J$. Let us assume that
$I_r\subset J$ and leave the other case to the reader.
For every $i\in I$,
we have $$i+\ln 2\sqrt{n}/\kappa\in I_r.$$ (To simplify notation,
we assume that $\ln2 \sqrt{n}/\kappa$ is a natural number.) 
Let $j:=i+\ln 2\sqrt{n}/\kappa$. By the assumption \eqref{Wkappa1},
we find that
\begin{equation}
\PP(W=j)\geq \PP(W=i)\left( 1-\frac{\kappa}{\sqrt{n}}  \right)^{\ln 2\sqrt{n}/\kappa} .
\end{equation}
The last inequality above yields
$$
\PP(W\in I_r)\geq \PP(W\in I)
\left( 1-\frac{\kappa}{\sqrt{n}}  \right)^{\ln 2\sqrt{n}/\kappa},
$$
and so
\begin{equation}
\label{IJ}
\PP(W\notin I)\geq \PP(W\in I)
\left( 1-\frac{\kappa}{\sqrt{n}}  \right)^{\ln 2\sqrt{n}/\kappa}.
\end{equation}
Note that
$$\lim_{n\rightarrow\infty}
\left( 1-\frac{\kappa}{\sqrt{n}}  \right)^{\ln 2\sqrt{n}/\kappa}
=\frac12,$$
and thus for $n$ large enough,
$$
\left( 1-\frac{\kappa}{\sqrt{n}}  \right)^{\ln 2\sqrt{n}/\kappa}
>\frac13.$$
The last inequality together with \eqref{IJ}, yields
$$\PP(W\notin I)\geq \frac{\PP(W\in I)}{3},$$ 
from which it follows that
\begin{equation}
\label{notin}
\PP(W\notin I)\geq \frac14.
\end{equation}
We have that
\begin{equation}
\label{VarWE2}
{\rm Var}\, W=\sum_{i\in J}p_i(i-\E\, W)^2
\geq \sum_{i\in J\backslash I}p_i(i-\E\, W)^2,
\end{equation}
where $p_i:=\PP(W=i)$.
Note that for every $i\notin I$, 
$$(i-\E\, W)^2\geq\left( \frac{\sqrt{n}\ln 2}{2\kappa}\right)^2.$$
Hence,
the right-hand side of  \eqref{VarWE2} is larger than
$$\PP(W\notin I)
n\frac{(\ln2)^2}{4\kappa^2}
.$$
This last fact combined
 with \eqref{notin} in \eqref{VarWE2}, leads to
$${\rm Var }W\geq n \frac{(\ln2)^2}{16\kappa^2}.
\eqno{\rule{0.5em}{0.5em}}$$
\end{eqproof}

The next lemma will also be useful:
\begin{lemma}
\label{3.4}There exists a constant
$c_M>0$, not depending on $\vec{M}$ or $n$ such that
\begin{equation}
\label{cM}
{\rm Var }(M_1\mid \vec{M},\vec{N}\in I^n)\geq c_Mn
\end{equation}
for all $\vec{M}$, such that
$\PP(\vec{M},\vec{N}\in I^n)>0$.
\end{lemma}
\begin{eqproof}Assume that $\vec{n}=(n_1,n_2,n_4,n_5)\in I^n$.
Let $m_1\in\mathbb{N}$ and $\vec{m}\in H_a\cup H_b$, be such that
$$\vec{n}=\vec{m}+m_1\vec{e}.$$
(Note that the above equality uniquely determines $\vec{m}$
 and
$m_1$.)
Let $I_1$ be the integer interval such that
$$\left(\mathbb{N}\cdot\vec{e}+\vec{m}\right)\cap I^n
=I_1\cdot\vec{e}+\vec{m}.$$
We also write $I(\vec{m})$ to indicate that 
$I_1$ only depends on $\vec{m}$. Note that the condition
$$\vec{M}=\vec{m},\quad \vec{N}\in I^n,$$
can now be rewritten as
$$\vec{M}=\vec{m},\quad M_1\in I_1.$$
The equation \eqref{y2} can be rewritten as:
\begin{equation}
\frac{\PP(\vec{M}=\vec{m},M_1=m_1+1)}
{\PP(\vec{M}=\vec{m},M_1=m_1)}=
\frac{(n_2+1)(n_4+1)}{n_1n_5}\,,
\end{equation}
and so
\begin{equation}
\label{MP}
\frac{\PP(M_1=m_1+1\mid \vec{M}=\vec{m},M_1\in I_1)}
{\PP(M_1=m_1\mid \vec{M}=\vec{m},M_1\in I_1)}=
\frac{(n_2+1)(n_4+1)}{n_1n_5}\,.
\end{equation}
We have that there exist $a,b,c,d>0$ not depending on $n$,
such that the following four inequalities hold
$$|\mu_1-r_1|<a\sqrt{n},$$
$$|\mu_2-r_2|<b\sqrt{n},$$
$$|\mu_4-r_4|<c\sqrt{n},$$
$$|\mu_5-r_5|<d\sqrt{n},$$
for any $\vec{r}=(r_1,r_2,r_4,r_5)\in I^n$.
Hence there exist $a^*,b^*,c^*,d^*>0$, such that
\begin{equation}
\label{rinequ}
\frac{(n_2+1)(n_4+1)}{n_1n_5}=
\frac{(\mu_2+b^*\sqrt{n})(\mu_4+c^*\sqrt{n})}
{(\mu_1+a^*\sqrt{n})(\mu_5+d^*\sqrt{n})},
\end{equation}
and such that
$|a^*|<a$, $|b^*|<b$, $|c^*|<c$ and $|d^*|<d.$
Note that
$$\mu_i=\frac{n}{2^i},$$
which implies that
\begin{equation}
\label{rinequ2}
\frac{(\mu_2+b^*\sqrt{n})(\mu_4+c^*\sqrt{n})}
{(\mu_1+a^*\sqrt{n})(\mu_5+d^*\sqrt{n})}
=\frac{(1/2^2+b^*/\sqrt{n})(1/2^4+c^*/\sqrt{n})}
{(1/2+a^*/\sqrt{n})(1/2^5+d^*/\sqrt{n})}\,.
\end{equation}

Let
$$f:(x,y,z,w)\mapsto\frac{(1/2^2+x)(1/2^4+y)}
{(1/2+z)(1/2^4+w)}.$$
Note that the maps $f$
and $1/f$ are both continuously differentiable
in an open neighborhood of $\vec{0}=(0,0,0,0)$ and that
$f(0,0,0, 0)=1$.
Hence, there exist a constant
$\kappa>0$ independent of $n$, such that for 
all $n$ large enough,  
\begin{equation}
\label{kappa*}
\frac{1}{1-\kappa\sqrt{n}}\geq\frac{(1/2^2+b^*/\sqrt{n})(1/2^4+c^*/\sqrt{n})}
{(1/2+a^*/\sqrt{n})(1/2^5+d^*/\sqrt{n})}\geq 1-\kappa\sqrt{n},
\end{equation}
which holds for every $a^*,b^*,c^*,d^*$ such that
$|a^*|\leq a$, $|b^*|\leq b$, $|c^*|\leq c$ and $|d^*|\leq d$.
Combining, \eqref{MP}, \eqref{rinequ}, \eqref{rinequ2} and \eqref{kappa*},
we find that
$$\frac{\PP(M_1=m_1+1\mid \vec{M}=\vec{m},M_1\in I_1)}
{\PP(M_1=m_1\mid \vec{M}=\vec{m},M_1\in I_1)}
\geq 1-\kappa \sqrt{n},
$$
and 
$$\frac{\PP(M_1=m_1\mid \vec{M}=\vec{m},M_1\in I_1)}
{\PP(M_1=m_1+1\mid \vec{M}=\vec{m},M_1\in I_1)}
\geq 1-\kappa \sqrt{n}.
$$
Note that the last two inequalities above are nothing but the conditions
\eqref{Wkappa} and \eqref{Wkappa1} in Lemma~\ref{3.3}. For this we take
for the random variable $W$ of Lemma~\ref{3.3}, the random variable $M_1$
conditional on $\vec{M}=\vec{m}$ and $M_1\in I_1$.
In order to apply Lemma~\ref{3.3}, we also need
to verify that the interval $I(\vec{m})$ has diameter of
at least $3\ln2\sqrt{n}/\kappa$. For this, note that
by choosing the constant $C>0$ in the Definition~\ref{defIn} 
of $I^n$ large enough (but not depending on
$n$), we get that
the diameter of $I(\vec{m})$ is larger than $3\ln2\sqrt{n}/\kappa$
for every $\vec{m}$, such that
 $$\PP(\vec{M}=\vec{m},\vec{N}\in I^n)>0.$$
Applying Lemma~\ref{3.3},
and obtain that
$${\rm Var}(M_1\mid \vec{M}=\vec{m},M_1\in I_1)=
{\rm Var}(M_1\mid \vec{M}=\vec{m},\vec{N}\in I^n)
\geq n\frac{(\ln 2)^2}{16\kappa^2}.\eqno{\rule{0.5em}{0.5em}}
 $$
\end{eqproof}

Next, applying \eqref{cM} into \eqref{fundi2}, we find that
conditional on $\vec{N}\in I^n$ and when $E^n_{\rm slope}$ holds,
\begin{equation}\label{fundi2b}
{\rm Var}(L(M_1,\vec{M})\mid L(\cdot ),\vec{M})\geq n\left(
10^{-4}\left(1-\frac{2n^{1/10}}{10^{-2}\sqrt{c_M n}}
\right)c_M\right).
\end{equation}
Using \eqref{fundi2b}, we see that
\eqref{alin1} becomes:
\begin{align}\label{alin3}
 \E&\left({\rm Var}(L(M_1,\vec{M})\mid L(\cdot ),\vec{M})\right)\\
&\label{alin4}\quad \geq n\left(
c_M 10^{-4}\left(1-\frac{2}{10^{-2}\sqrt{c_M} n^{4/10}}
\right)\right)
\PP(\vec{N}\in I^n) \PP(E^n_{\rm slope}).
\end{align}
With the help of \eqref{beginning} and of
\eqref{Ln=LN}, we then find
\begin{equation}\label{final}
{\rm Var }L_n\geq n\left(
c_M 10^{-4}\left(1-\frac{2}{10^{-2}\sqrt{c_M} n^{4/10}}
\right)\right)
 \PP(\vec{N}\in I^n) \PP(E^n_{\rm slope}).
\end{equation}
To finish the proof of Theorem~\ref{mainresult} it suffices
to prove that $\PP(\vec{N}\in I^n)$ and $\PP(E^n_{\rm slope})$
are both uniformly bounded below by a constant as $n$ goes to infinity.
This is the content of the two next lemmas.
\begin{lemma}
We have 
$$\PP(\vec{N}\in I^n)\geq \frac15.$$
\end{lemma}
\begin{proof}Let $Z_i$ be the indicator function which is
equal to one if 
$$X_i=1,X_{i+1}=0,X_{i+2}=1.$$
and $Z_i=0$ otherwise.  In other words, $X_i=1$ if there is a block
of one zero in $X$ starting at the position $i$.  Then,
$$N_1=\sum_{i=1}^{n-3}Z_i,$$
and hence
\begin{equation}
\label{VarN1}
{\rm Var }N_1=\E\left(\sum_{i=1}^{n-3}(Z_i-\E\, Z_i)\right)^2=
\sum_{i,j}\E\left((Z_i-\E\, Z_i)(Z_j-\E\, Z_j)\right).
\end{equation}
When $|i-j|>2$, then $Z_i$ and $Z_j$ are independent of each other
and so
$$\E\left((Z_i-\E\, Z_i)(Z_j-\E\, Z_j)\right)=0.$$
This implies that in the sum on the very right side of
\eqref{VarN1}, there at most $5n$ terms which are different from
zero.  By the Cauchy-Schwarz inequality, 
$$
\E\left((Z_i-\E\, Z_i)(Z_j-\E\, Z_j)\right)\leq
{\rm Var}\, Z_1,
$$
hence, for all $i,j\in[0,n-3]$,
\begin{equation}
\label{Cauchyschwarz}
\E\left((Z_i-\E\, Z_i)(Z_j-\E\, Z_j)\right)\leq
\frac14 .
\end{equation}
Since there are less than $5n$ non-zero terms in the sum on the right
side of \eqref{VarN1}, \eqref{Cauchyschwarz} leads to:
\begin{equation}
\label{boundVarN1}
{\rm Var }N_1\leq \frac{5n}{4}\,.
\end{equation}
Similarly,
\begin{equation}
\label{boundsVarN1234}
{\rm Var }N_2\leq \frac{5n}{4},\quad
{\rm Var }N_4\leq \frac{9n}{4},\quad
{\rm Var }N_5\leq \frac{11n}{4}\,.
\end{equation}
In \eqref{IC}, we have chosen $I^n$ such that
$$\bigcap_{i=1,2,4,5}\left\{N_i\in\left[\mu_i
-\sqrt{55n/4},\mu_i+\sqrt{55n/4}\right]\right\}
\subset\left\{\vec{N}\in I^n \right\},$$
and hence
\begin{equation}\label{probsnotin}
\PP(\vec{N}\notin I^n)\leq
\sum_{i=1,2,4,5}\PP
\left(N_i\notin\left[\mu_i-\sqrt{55n/4},\mu_i+\sqrt{55n/4}\right]\right).
\end{equation}
The Bienaym\'e-Chebycheff inequality yields that
$$
\PP \left(N_i
\left[\mu_i-\sqrt{55n/4},\mu_i+\sqrt{55n/4}\right]\right)\leq 
\frac{{\rm Var}\, N_i}{55n/4},$$
which together with \eqref{boundsVarN1234} and \eqref{probsnotin},
yield:
$$\PP(\vec{N}\notin I^n)\leq
\sum_{i=1,2,4,5}\left(\frac{11n}{4}\right)\left(\frac{4}{55n}\right)
=\frac45\,,$$
finishing the proof.
\end{proof}

In the next lemma,
recall that $\mu_i$ is the expected number of blocks
of zeros of length $i$ in $X$.
\begin{lemma}
For all $n$ and all $i$, we have that
\begin{equation}
\label{muiless}
\left| \mu_i-n\left(\frac12   \right)^{i+2}\right|\leq i.
\end{equation}
\end{lemma}
\begin{proof}
Let $R_j^i$ be the Bernoulli random variable which indicates
if there is  a block of zeros of length $i$
starting at the point $j$ in $X$.  In other words,
$R_j^i=1$ if
$$X_j=X_{j+1}=\cdots=X_{j+i-1}=0\;\;{\rm and}\;\;
X_{j-1}=X_{j+i}=1$$
 for $j\in[2,n-i]$.
For $j=1$ we have  $R_j^i=1$ if
$$X_1=X_{2}=\cdots=X_{i}=0\;\;{\rm and}\;\;
X_{i+1}=1.$$
For $j=n-i+1$, we have that $R_j^i=1$ if
$$X_{n-i+1}=X_{n-i+2}=\cdots=X_{n}=0\;\;{\rm and}\;\;
X_{n-i}=1.$$
We find that
$$
N_i=\sum_{j=1}^{n-i+1}R_j^i,
$$
and hence
\begin{equation}
\label{N}
\mu_i=\E\, N_i=\sum_{j=1}^{n-i+1}\E\, R_j^i.
\end{equation}
Note that for $j\in[2,n-i]$, we have that
$\E\, R_j^i=(1/2)^{i+2}$.
We also have that $\E\, R_1^i=\E\, R_{n-i+1}^i=(1/2)^{i+1}$.
This with \eqref{N}, yields
\begin{equation}
\label{Nin}
\mu_i=(n-i)(1/2)^{i+2}+2(1/2)^{i+1}.
\end{equation}
This last equality directly implies \eqref{muiless}.
\end{proof}

Our next lemma is related to the
local limit theorem for multidimensional 
regenerative processes:
\begin{lemma}
There exist a constant $k_1>0$ such that for 
all $n$, and
all $\vec{n}\in I^n$, 
\begin{equation}\label{k_1}
\PP(\vec{N}=\vec{n})\geq n^{-k_1}.
\end{equation}
\end{lemma} 
\begin{proof}Let $T_i$ be the location of the end of the $i$-th block in
the infinite sequence $X_1,X_2,X_3,\ldots$. 
To simplify notation, we assume that the finite sequence $X=X_1X_2\cdots X_n$
got extended to an infinite sequence 
of iid Bernoulli random variables with parameter $1/2$.  
In this way, $T_i$ is well defined even
when $T_i>n$. (Note that
here we consider all the blocks, and not just the blocks of zeros.)
Let $Z_i$ denote the length of the $i$-th block of $X$.  Hence, 
$Z_i=T_i-T_{i-1}$ and
$$T_i=Z_1+Z_2+\cdots+Z_i.$$ 

Note that $T_1,T_2,\ldots$ are the arrival times
of a renewal process.  The interarrival times $Z_1,Z_2,\ldots$
are iid geometric random variables with parameter $2$ and hence $\E[Z_i]=2$.
We also assume that $n/4$ is an integer in order to simplify notation.
There are $n/2$ blocks in $X$ if and only if $T_{n/2}=n$.
The event that the blocks of ones cover half the text $X$ can be described
by the equality
\begin{equation}
\label{half}
\sum_{i=1}^{n}X_i=\frac{n}{2}.
\end{equation}
Let $\vec{n}=(n_1,n_2,n_4,n_5)\in I^n$.  We have that
\begin{equation}
\label{vecn}
\PP(\vec{N}=\vec{n})\geq \PP\left(\vec{N}=\vec{n},\sum_{i=1}^{n}X_i=\frac{n}{2},
T_{n/2}=n\right).
\end{equation}
Let $V^j_i$ be the indicator variable which is equal to one
if the $i$-th block of zeros in $X$ had length $j$
 and $V^j_i=0$ otherwise.  Let
$$\vec{V}_i:=(V^1_i,V^2_i,V^4_i,V^5_i).$$
Note that $\vec{V}_1,\vec{V}_2,\ldots$ form a sequence
of iid random vectors.
The event
$$\left\{\vec{N}=\vec{n},\sum_{i=1}^{n}X_i=\frac{n}{2},T_{n/2}=n\right\}$$
can be interpreted as the event that there are $n/2$ blocks in
$X$ and that the blocks of ones cover half the text and
 that $\vec{N}=\vec{n}$. Its probability can be computed
as follows: First compute the contribution due to the fact that
the blocks of ones cover half the text.  Half of all blocks
are blocks of ones.  Hence we have exactly $n/4$ blocks of ones.
These blocks are iid with geometric distribution with parameter
$1/2$.
This gives us a factor
\begin{equation}
\label{n/4}
\PP(Z_1+Z_2+\cdots+Z_{n/4}=n/2)
\end{equation}
Second, we compute the probability that
among the first $n/4$ blocks of zeros we have the right number of blocks
of length $1$, $2$, $4$ and $5$.  We find
\begin{equation}
\label{probV}
\PP(\vec{V}_1+\cdots+\vec{V}_{n/4}=\vec{n}).
\end{equation}
Finally we compute the probability that the remaining $n/4-n_1-n_2-n_4-n_5$
blocks of zeros, cover up a total length of $n/2-n_1-2n_2-4n_4-5n_5$. 
This yields the probability
\begin{equation}\label{probW}
\PP(W_1+W_2+\cdots+W_{n^*}=n/2-n_1-2n_2-4n_4-5n_5),
\end{equation}
where $n^*:=n/4-n_1-n_2-n_4-n_5$ and
$W_1,W_2,\ldots$ are iid random variables with distribution
$$\mathcal{L}(Z_i\mid Z_i\notin\{1,2,4,5\}).$$
Summarizing, we find
\begin{align*}
&\PP\left(\vec{N}=\vec{n},\sum_{i=1}^{n}X_i=\frac{n}{2},
T_{n/2}=n\right)=\\
&\PP(Z_1+\cdots+Z_{n/4}=n/2)
\PP(\vec{V}_1+\cdots+\vec{V}_{n/4}=\vec{n})
\PP(W_1+\cdots+W_{n^*}=n/2-\bar{n}),
\end{align*}
where $\bar{n}=n_1+2n_2+4n_4+5n_5$
and hence with the help of \eqref{vecn}:
\begin{equation*}
\PP(\vec{N}=\vec{n})\geq
\PP\left(Z_1+\cdots+Z_{n/4}=\frac n2\right)
\PP(\vec{V}_1+\cdots+\vec{V}_{n/4}=\vec{n})
\PP(W_1+\cdots+W_{n^*}=\frac n2-\bar{n}).
\end{equation*}
Since $(n/4)\E\, Z_1=n/2$,  by the 
local limit theorem, 
there exists $k_2>0$ not depending on $n$ such that
\begin{equation}
\label{k_2}
\PP(Z_1+\cdots+Z_{n/4}=n/2)\geq \frac{k_2}{\sqrt n}
\end{equation}
for all $n$.
Note that 
\begin{equation}
\label{Vi}
\E\, V_1^i=\left(\frac12  \right)^i,
\end{equation}
and so, together with the inequality \eqref{muiless},
this yields:
\begin{equation}
\label{muE}
|\vec{\mu}-(n/4)\E\, \vec{V}_1|\leq 1+2+4+5=12,
\end{equation}
where $|\cdot|$ is the $\ell^1$-norm
in $\mathbb{R}^4$ for which $|(x,y,z,w)|=|x|+|y|+|z|+|w|$.
Then,
\begin{equation}
\label{+}
|\vec{n}-(n/4)\E\, \vec{V}_1|
\leq |\vec{n}-\vec{\mu}|+|\vec{\mu}-(n/4)\E\, \vec{V}_1|,
\end{equation}
and since $\vec{n}\in I^n$, there exists a
constant $K>0$ independent of $n$ or $\vec{n}\in I^n$
such that
\begin{equation}\label{Ksqrt}
|\vec{n}-\vec{\mu}|\leq K \sqrt{n}.
\end{equation} 
Combining \eqref{muE}, \eqref{+} and \eqref{Ksqrt} yields
$$|\vec{n}-(n/4)\E\, \vec{V}_1|\leq K\sqrt{n}+12
$$
which hold for every $n$ and every $\vec{n}\in I^n$.
 From the last inequality above, with the help of the 
local limit theorem, we obtain that there exists
$k_3>0$ not depending on $n$ or $\vec{n}\in I^n$ such that
\begin{equation}
\label{k_3}
\PP(\vec{V}_1+\cdots+\vec{V}_{n/4}=\vec{n})\geq \frac{k_3}{\sqrt n}.
\end{equation}
Next, we want to prove that
\begin{equation}
\label{lastlast}
\left|n^*\E\, W_1-n/2+n_1+n_2+n_4+n_5   \right|
\end{equation}
is also of order $\sqrt{n}$.
Let $$u_i:=\frac{n}{4}\left(\frac12\right)^i.$$
We already know that for a constant $K>0$, independent of
$n$ or $\vec{n}\in I^n$, 
$$|\mu-\vec{n}|\leq K\sqrt{n}.$$
Hence using also \eqref{muiless}, it follows that there exists
$K_2>0$ such that the difference between \eqref{lastlast}
and
\begin{equation}
\label{lastlast2}
\left|(n/4-u_1-u_2-u_4-u_5)\E\, W_1-n/2+u_1+u_2+u_4+u_5\right|
\end{equation}  
is less than $K_2\sqrt{n}$.
Now consider a sequence of $n/4$ iid geometric
random variables with parameter $1/2$. The expectation
of the sum of these random variables is:
$$\E\, T_{n/4}=\frac{n}{4}\E\, Z_1=\frac{n}{2}.$$
Let $W$ denote the sub-sum obtained by only taking the terms 
not equal to 1, 2, 4 or 5. Let $\widetilde{W}$ denote the sub-sum obtained
by taking only those terms in the sum which are equal
to 1, 2, 4 or 5.
Let $b$ denote the total number of random variables, among
our collection of $n/4$, that take a value equal
to 1, 2, 4 or 5. 
Note that the probability for any one of such variables to take on the value
$i$ is equal to $(1/2)^i$.  Hence, the expected number of variables
among our set of $n/4$ which take on the value $i$ is equal to 
$(n/4)(1/2)^i=u_i$.
We find
$$\E\, T_{n/4}= \E(W+\widetilde{W})=
\E\, W+\E\widetilde{W},$$
but
$$\E\, W=\E\, W_1 \E((n/4)-b)=\E\, W_1((n/4)-u_1-u_2-u_4-u_5)$$
and
$$\E\, \widetilde{W}=u_1+2u_2+4u_4+5u_5.$$
Combining the last three inequalities yields:
$$
n/2=\E\, T_{n/4}=
\E\, W_1((n/4)-u_1-u_2-u_4-u_5)+u_1+2u_2+4u_4+5u_5,
$$
and thus
\begin{equation}
\label{tired}
0=\E\, W_1((n/4)-u_1-u_2-u_4-u_5)-(n/2)+u_1+2u_2+4u_4+5u_5.
\end{equation}
Therefore, the expression \eqref{lastlast2} is equal to zero and so
$$\left|n^*\E\, W_1-n/2+n_1+n_2+n_4+n_5   \right|\leq K_2\sqrt{n}
.$$
The above inequality combined with the local limit theorem
yields that there exists $k_4>0$, not depending on $n$ or
$\vec{n}\in I^n$ such that
\begin{equation}
\label{k_4}
\PP(W_1+W_2+\cdots+W_{n^*}=n/2-n_1-2n_2-4n_4-5n_5)\geq k_4 \sqrt n
\end{equation}
The inequalities \eqref{k_2}, \eqref{k_3} and \eqref{k_4} 
together imply \eqref{k_1}.
\end{proof}

The next lemma is proved assuming that Theorem~\ref{maintheq}
holds. (In turn, Theorem~\ref{maintheq} is proved in 
Subsection~\ref{overview}.)

\begin{lemma}There exists $s^*>0$ independent of
$n$ such that if $s(1,1)\geq s^*$, then
$$\PP(E^n_{\rm slope})\rightarrow 1,$$
as $n\rightarrow\infty$.
\end{lemma}
\begin{proof}
Let $\vec{m}\in H_a\cup H_b$, be such that
\begin{equation}
\label{conditionvecm}
\PP(\vec{M}=\vec{m},\vec{N}\in I^n)>0.
\end{equation}
Let $E^n_{\vec{m}}$ be the event that
$E^n_{\rm slope}$ holds on the subset $I(\vec{m}) \vec{e}+\vec{m}$.
(The set $I(\vec{m})$ got defined in the proof of Lemma~\ref{3.4}.)
More precisely,  $E^n_{\vec{m}}$ 
is the event that
for all $i,j\in I(\vec{m})$, such that 
$j-i\geq n^{1/10}$,
$$L(\vec{m}+j\vec{e})-L(\vec{m}+i\vec{e})\geq 10^{-2} |j-i|.$$
We find that
\begin{equation}
\label{Eslopem}
E^n_{\rm slope}=\bigcap_{\vec{m}}E^n_{\vec{m}},
\end{equation}
where the intersection is taken over all
$\vec{m}\in H_a\cup H_b$, such that
\eqref{conditionvecm} holds. Note that there exists
a constant $c_I$ (not depending on $n$)
such that for all $n$ there are less than $c_I n^2$
points in the set $I^n$.  Hence there are also less
than $c_I n^2$ vectors
 $\vec{m}\in H_a\cup H_b$ satisfying
\eqref{conditionvecm}. It follows that in the intersection on the right side
of \eqref{Eslopem}, there are less than $c_I n^2$ terms.
Now, \eqref{Eslopem} implies that
\begin{equation}
\label{Eslopem2}
\PP((E^{n}_{\rm slope})^c)\leq \sum_{\vec{m}}\PP((E^{n}_{\vec{m}})^c),
\end{equation}
where the sum, above, is taken over
all 
$\vec{m}\in H_a\cup H_b$, such that
\eqref{conditionvecm} holds.  There are less than
$c_I n^2$ terms in the sum on the right side of \eqref{Eslopem2}.
Hence to prove the lemma, we only need an exponentially small upper-bound 
independent of $\vec m$
for $\PP((E^{n}_{\vec{m}})^c)$.
Let $A^n_i$ denote the event ``that $A^n$ holds for $X(\vec{m}+i\vec{e})$."
More precisely, $A^n_i$ is the event that the following two conditions
hold:
\begin{equation}
\label{probcondA1m}
\PP(L(\vec{m}+(i+1)\vec{e})-L(\vec{m}+i\vec{e})= 1\mid X(\vec{m}+i\vec{e}),Y)
\geq \frac{31}{128}
- \epsilon_1,
\end{equation}
\begin{equation}
\label{probcondA3m}
\PP(L(\vec{m}+(i+1)\vec{e})-L(\vec{m}+i\vec{e})=- 1\mid X(\vec{m}+i\vec{e}),Y)
\leq\frac{1}{32}
+\epsilon_1.
\end{equation}
Note that the difference between $L(\vec{m}+(i+1)\vec{e})$
and $L(\vec{m}+(i+1)\vec{e})$ is at most one. We assume that
the constant $\epsilon_1\leq 7/32$. The inequalities
\eqref{probcondA1m} and \eqref{probcondA3m} then give that
for every $(x,y)\in A^n_i$, we have
\begin{equation}
\label{biasAni}
\E\left(L(\vec{m}+(i+1)\vec{e})-L(\vec{m}+i\vec{e})= 1
\mid X(\vec{m}+i\vec{e})=x,Y=y\right)
\geq \frac12\,.
\end{equation}
From a positive bias like in 
\eqref{biasAni}, one can hope to prove that
the event $E^n_{\vec{m}}$ holds with probability one minus
an exponentially small quantity. The only
problem is that inequality \eqref{biasAni} holds only for
$(x,y)\in A^n_i$.  If we condition on $\cap_{i\in I(\vec{m})}A^n_i$,
 we introduce complicated dependencies so that
we can no longer use large deviation results for martingales.  The trick
is to introduce help-variables $Y_i$.  When,
$A^n_i$ holds let 
$$Y_i:=L(\vec{m}+(i+1)\vec{e})-L(\vec{m}+i\vec{e}),$$
otherwise let $Y_i:=1$.  We have that
$Y_i$ is $\sigma(Y,X(\vec{m}+j\vec{e})|j\leq i)$-measurable.
Also, because of \eqref{biasAni}, we have that almost surely,
\begin{equation}
\label{biasY}
\E\left(Y_i\mid Y,X(\vec{m}+j\vec{e})\right)\geq\frac12.
\end{equation}
Let $E^n_{\vec{m},Y}$ be the event that
for all $i,j\in I(\vec{m})$, such that 
$j-i\geq n^{1/10}$,
we have that
$$\sum_{k=i}^{j}Y_k\geq 10^{-2} |j-i|.$$
Note that when $\bigcap_{i\in I(\vec{m})}A^n_i$ holds,
then the events $E^n_{\vec{m},Y}$ and $E^n_{\vec{m}}$ are
identical. Hence,
$$
 \bigcap_{i\in I(\vec{m})}A^n_i \;\cap\;E^n_{\vec{m},Y}
\subset E^n_{\vec{m}},
$$
and thus
\begin{equation}
\label{Ani}
\PP((E^{n}_{\vec{m}})^c)\leq \PP((E^{n}_{\vec{m},Y})^c)+
\sum_{i\in I(\vec{m})}\PP((A^{n}_i)^c).
\end{equation}
Note that $X(\vec{m}+i\vec{e})$ has distribution
$\mathcal{L}(
X\mid \vec{N}=\vec{m}+i\vec{e})$.  This implies that
$$\PP((A^{n}_i)^c)=
\frac{\PP\left((A^{n})^c\cap \left\{\vec{N}=\vec{m}+i\vec{e}\right\}\right)}
{\PP(\vec{N}=\vec{m}+i\vec{e})}
\leq
\frac{\PP((A^{n})^c)}{\PP(\vec{N}=\vec{m}+i\vec{e})}\,.
$$
Using the last inequality with \eqref{k_1}, we obtain
$$\PP((A^{n}_i)^c)\leq
\PP((A^{n})^c) n^{k_1}.$$
The last inequality implies:
\begin{equation}
\label{avant}
\sum_{i\in I(\vec{m})}\PP((A^{n}_i)^c)\leq c_I n^{2+k_1}\PP((A^{n})^c).
\end{equation}
Classical exponential inequalities and (\ref{biasY}) 
show that
\begin{equation}
\PP((E^{n}_{\vec{m},Y})^c)\leq e^{-k_Y\cdot n^{1/10}},
\end{equation}
where $k_Y>0$ is a constant independent of $n$.
Next, \eqref{Ani} and \eqref{avant} together imply that
\begin{equation}
\label{Ani2}
\PP((E^{n}_{\vec{m}})^c)\leq e^{-k_Y n^{1/10}}+
c_I n^{2+k_1}\PP((A^{n})^c).
\end{equation}
We can now plug \eqref{final} into \eqref{Ani2} to obtain
$$\PP((E^{n}_{\vec{m}})^c)\leq e^{-k_Y n^{1/10}}+
c_I n^{2+k_1}e^{-c_1n}.
$$
The upper-bound in the last inequality above is
exponentially small in $n^{1/10}$, and so by 
\eqref{Eslopem2}, 
$\PP((E^{n}_{\rm slope})^c)$ must also be exponentially small
in $n^{1/10}$.
\end{proof}

\section{Combinatorics}
We already mentioned that
$\PP(X_i=1)=1/5$. Let $\epsilon>0$ be a small quantity not
depending on $n$.
Let $B^n_0$ be the event
that in both $X$ and  $Y$ there are about $n/2$ ones.
  More precisely,  $B^n_0$ is the event that
$$\left|\sum_{i=1}^{n}X_i-\frac n2\right|\leq \frac{\epsilon n}{16}$$
and 
$$\left|\sum_{i=1}^{n}Y_i-\frac n2\right|\leq \frac{\epsilon n}{16}$$
both hold.
Let $B^n_1$ be the event that  any optimal alignment
of $X$ and $Y$ contains at least $ n/2-\epsilon n/8$ pairs of aligned
ones.
\begin{lemma}\label{lem4.1}Let $\epsilon>0$, and let  
\begin{equation}
\label{conditionepsi0}
s(1,1)\ge \epsilon/16,
\end{equation}
then, for any $n$,
\begin{equation}
\label{B0B1}
B^n_0\subset B^n_1.
\end{equation}
\end{lemma}
\begin{proof}
Let $v_1$ be the alignment of $X$ and $Y$
which aligns only ones and as many as possible.
Let $S_1$ denote the score obtained by aligning
$X$ with $Y$ via $v_1$. When, $B^n_0$ holds, then
$$
S_1\geq s(1,1) \left(\frac n2-\frac{\epsilon n}{16}\right),
$$
and therefore
\begin{equation}
\label{s1}
L_n\geq s(1,1)\left(\frac n2-\frac{\epsilon n}{16}\right)>
\left(\frac n2-\frac{\epsilon n}8\right) s(1,1)+n,
\end{equation}
since
$s(1,1)>16/\epsilon$.

Recall that $s(0,0)=1$.  Hence, if the texts $X$ and $Y$
consist only of zeros the maximum score would be equal to $n$.
This also implies that $n$ is an upper bound for the contribution
made by the aligned zeros to the score of any alignment.
Assume now that $v$ is an alignment which aligns no more that
$n/2-\epsilon n/8$, pairs of ones. Let $S_v$ denote the score
of $v$.  
When $B^n_0$ holds, using the bound for the contributions of the zeros
in the score, it follows that
\begin{equation}
\label{scoreless}
S_v\leq\left(\frac n2-\frac{\epsilon n}8\right) s(1,1)+n.
\end{equation}
Together, \eqref{s1} and \eqref{scoreless} imply that
$v$ is not an optimal alignment.  This implies that
when $B^n_0$ holds,
then any optimal alignment contains at least
$n/2-\epsilon n/8$ pairs of aligned ones.
In other words, 
$B^n_0$ implies $B^n_1$, when \eqref{conditionepsi0} holds.
\end{proof}

Recall now that the score $s(1,1)$ is taken to be large. This
 ensures, that typically there is a large 
proportion of the total number of ones, which get matched with a one
by the optimal alignment. 
We introduce a special notation for alignments, which is
convenient to describe alignments which align most ones with ones.
Let us start with a numerical example:

Take the finite sequence of pairs of natural numbers
$$(0,0),(0,1),(1,0).$$
According to our notations, this sequence
represents an alignment which does the following:
\begin{itemize}
\item First $(0,0)$ indicates that the first one of $X$ is aligned with
the first one of $Y$ without skipping any one.
\item The second pair $(0,1)$ indicates that after the first pair
of aligned ones, a one is skipped in the
$Y$-sequence and no one in the $X$-sequence.
\item The third pair $(1,0)$ indicates that after the second pair of aligned
ones, a one is skipped in $X$ and no one in $Y$.\end{itemize}

Take for example the sequences $X=0101011$ and $Y=0001111$. 
The alignment $v=((0,0),(0,1),(1,0))$ is then equal to:
$$\begin{array}{c|c|c|c|c|c|c|c|c}
0& & &1&0&1&0&1&1\\\hline
0&0&0&1&1&1& & &1
\end{array}\,.
$$
Recall that the score for aligning a zero with a one is zero:
$s(0,1)=0$.
Let $V^k$ denote the set of alignments which align exactly
 $k$ pairs 
of ones with each other and such that there is a proportion of less than $\epsilon/2$
ones not belonging to pairs of aligned ones. (The ones which are not aligned
with ones are counted up to the last pair of aligned ones.)
 Hence, with our representation
of alignments of pairs of ones as sequences of couples of natural numbers
we find
$$V^k:=\left\{(v_1,v_2,\ldots,v_k): v_1,v_2,\ldots,v_k\in
\mathbb{N}\times\mathbb{N},|v_1|+|v_2|+\cdots+|v_k|\leq \epsilon k/2\right\},
$$
where if $v=(a,b)\in\mathbb{N}^2$, 
$|v|:=a+b$.
Let 
$$V:=\bigcup_{k\geq p^*n}V^k$$
where $p^*:=1/2-\epsilon/8$.
Let $B^n_2$ be the event that  any optimal alignment of
$X$ with $Y$ is contained in $V$.  In other words, $B^n_2$
holds when for every alignment such that $S_v=L_n$,
we have $v\in V$. (Here $S_v$ is the score obtained
by aligning $X$ with $Y$ and using for this the alignment $v$.)

The next lemma shows that $B^n_0$ and $B^n_1$ together imply $B^2_n$:
\begin{lemma}Let $0<\epsilon<1$, then for any $n$,
\begin{equation}
\label{B1B2}
B^n_0\cap B^n_1\subset B^n_2.
\end{equation}
\end{lemma}
\begin{proof}Let $v$ be an alignment. 
We say that a one is {\it matched} by $v$ if $v$
aligns it with another one.  When it is clear from
the context which alignment $v$ we are talking about, we simply 
say that a one is matched. If $B^n_0$ and $B^n_1$ both hold, there are at most
$(\epsilon n/16)+(\epsilon n/8)$ non-matched ones in each word
$X$ and $Y$ for any optimal alignment $v$.
 On the other hand, $B^n_1$ ensures that at least
$n/2-\epsilon n/8$ ones are matched in each text with
a one from the other text by any optimal alignment $v$.  
This ensures that the proportion
of unmatched ones by the total number of matched ones is 
smaller or equal to
\begin{equation}
\label{fract}
\frac{\frac{\epsilon n}{16}+
\frac{\epsilon n}{8}}{\frac n2-\frac{\epsilon n}8}
=\frac{3\epsilon}{8}
\left(\frac{1}{1-\epsilon/4}\right).
\end{equation}
Hence by choosing  $\epsilon>0$  small enough, we get
\begin{equation}
\label{conditionepsi1}
\frac{3\epsilon}{8}
\left(\frac{1}{1-\epsilon/4}\right)\leq \frac{\epsilon}{2}.
\end{equation}
 We assume henceforth that the inequality \eqref{conditionepsi1}
holds, i.e., that $0<\epsilon <1$. 
This implies that the optimal alignment has a proportion
of non-matched ones to matched ones smaller than
$\epsilon/2$. Adding to this, that by $B^n_1$ there
are at least $n/2-\epsilon n/8$ pairs of aligned ones,
gives that any optimal alignment is in $V$.  Hence,
$B^n_2$ holds.  Therefore, $B^n_0$ and
$B^n_1$ together imply $B^n_2$, when $\epsilon>0$
is small enough.~\end{proof}

\medskip
For $v\in V^k$, let $(\pi_v(i),\nu_v(i))$ be the indices
of the $i$-th pair of ones aligned by $v$. Hence, if
$\pi_v(i)=j$ and $\nu_v(i)=k$, then all of the following properties hold:
\begin{itemize}
\item $X_j$ gets aligned with $Y_k$
\item $X_j=Y_k=1$
\item The $i$-pair of aligned ones (by $v$) is $X_j$ and $Y_k$.
\end{itemize}
For example, in the previous numerical example 
$\pi_v(1)=2$, $\pi_v(2)=4$ and $\pi_v(3)=7$. Furthermore,
$\nu_v(1)=4$, $\nu_v(2)=6$ and $\nu_v(3)=7$.

Let $v$ be an alignment. Recall that a one that gets aligned by
$v$ to another one (instead of aligned to a gap or a zero), is 
said to have been
{\it matched} by $v$.
Let $N_5(v)$ denote in the alignment $v$, the total number of
subsequent pairs of ones in $X$ satisfying all of the following conditions:
\begin{itemize}
\item The ones are both matched, that is aligned by $v$ with a one from $Y$.
\item Between the pair of ones in the $X$ text there are only zeros.
More precisely, we require that there is a block of five zeros, between
the pair of ones in the $X$-string.
\item The pair of ones in the $Y$-string with which our pair from $X$ is
aligned, should contain only zeros in between them.\end{itemize}
 More precisely: for
 $v\in V^k$, define:
$$N_5(v):=|\left\{
i<k:v_i=(0,0),\pi_v(i+1)-\pi_v(i)=6\right\}|,$$
where $v:=(v_1,v_2,\ldots,v_k)$.
Let $N_{5<}(v)$ denote in the alignment $v$, the total number of
subsequent pairs of ones in $X$ satisfying all of the following conditions:
\begin{itemize}
\item The ones are both matched, that is aligned by $v$ with a one from $Y$.
\item Between the pair of ones in the $X$ text there are only zeros
and moreover exactly five zeros.
\item The pair of ones in the $Y$ text with which our pair from $X$ is
aligned, should contain only zeros in between them and
 contain  strictly less than
$5$.\end{itemize}  
More precisely:
$$N_{5<}(v):=|\left\{
i<k:v_i=(0,0),\pi(i+1)-\pi(i)=6,\nu_v(i+1)-\nu_v(i)<6\right\}|.$$
Let $C^n$ be the event that
for all $v\in V$ which is an optimal alignment, 
$$\frac{N_{5<}(v)}{N_5(v)}\geq\frac{31}{32}- \frac{\epsilon_1}4\,.$$
Let $B^n_3$ be the event that in the sequence $X$ there are
at least $(1/32-\epsilon/16)n$ blocks of zeros of length
five.
Let $p_5(v)$ be the conditional probability on $X$, that
when picking a block of five zeros  at random in $X$,
 this block
happens to satisfy   the following two conditions:\begin{itemize}
\item[1)] The block is contained between two consecutive matched ones.
(Matched by the alignment $v$.)
In other words, between the two matched ones there 
is the block of length five and nothing else. 
\item[2)] The pair of ones in the $Y$-text to which the pair of consecutive ones
are aligned by $v$ contains only zeros in between them and strictly less than
five of them.\end{itemize}
In other words, $p_5(v)$ is the
conditional probability (conditional on $X$) that when picking
at random a block of zeros of length five
in $X$,  there exists
$i\leq k$, such that the randomly selected block
is equal to $[\pi_v(i-1)+1,\pi_v(i)-1]$ and all of the following properties hold:
$$\pi_v(i)-\pi_v(i-1)=6,$$
$$\eta_v(i)-\eta_v(i-1)<6,$$
and $|v_i|=0$,
where $$v=(v_1,v_2,\ldots,v_k)\in(\mathbb{N}\times\mathbb{N})^k.$$
 Recall that the total number of blocks with five zeros in $X$ is
denoted by $n_5$. Furthermore, each block of length five is selected
with equal probability among all blocks of five zero in $X$.
Hence, each block of five zeros in $X$ has a conditional
probability of $1/n_5$ to get selected.  There are
$N_{5<}$ blocks of five zeros in $X$ satisfying the conditions,
hence the conditional probability $p_5(v)$ is equal to: 
$$p_5(v):=\frac{N_{5<}(v)}{n_5}.$$

Let us return to a numerical example.
Let $X=1010101000001$ and let
$Y=1010001100011$.  Let $v$ be the alignment
$$
\begin{array}{c|c|c|c|c|c|c|c|c|c|c|c|c|c|c|c|c|c}
1&0&1&0& & &1&0&1&1&0&0&0&0&0&1&\\\hline
1&0&1&0&0&0&1& & &1&0&0&0& & &1 &1
\end{array}\,.
$$
In $X$ there is one block of five zeros.
Hence, $n_5=1$.
Moreover, this block is contained directly between
matched ones and the corresponding ones in the $Y$-text only
contain between them  zeros.  Hence, this block of
five zeros counts towards $N_5(v)$. We have
$N_5(v)=1$.
 The block  of five zeros in $X$ is matched with a block with
three zeros in $Y$, hence with a block having a strictly smaller number of zeros.
Thus, the block with five zeros is also counted towards
$N_{5<}$ and thus $N_{5<}=1$.  Hence, the conditional probability
$p_5$ is equal to $1$.

Let $E^n$ be the event
that
$$p_5(v)\geq\frac{31}{32}-\frac{\epsilon_1}{2}
$$ for every optimal alignment $v\in V$.
\begin{lemma}
Let $\epsilon>0$ and $\epsilon_1>0$ 
satisfy the inequality \eqref{conditionepsi2} below,
then for any $n$,
\begin{equation}
\label{B1En}
B^n_0\cap B^n_1\cap B^n_3\cap C^n\subset E^n.
\end{equation}
\end{lemma}
\begin{proof}
For any alignment $v\in V$, 
\begin{equation}
\label{williams}
\frac{N_{5<}(v)}{n_5}=\frac{N_{5<}(v)}{N_5(v)}
\frac{N_5(v)}{n_5}.
\end{equation}
When, the event $C^n$ holds, then for every 
$v\in V$, 
$$\frac{N_{5<}(v)}{N_5(v)}\geq\frac{31}{32}- \frac{\epsilon_1}4.$$
Combining the above inequality with \eqref{williams}
yields
\begin{equation}
\label{williamsII}
\frac{N_{5<}(v)}{n_5}\geq\left(
\frac{31}{32}- \frac{\epsilon_1}4
\right)
\frac{N_5(v)}{n_5}.
\end{equation}
Now, $B^n_0$ and $B^n_1$ together imply $B^n_2$,
hence when $B^n_0$ and $B^n_1$ both hold, then every
optimal alignment is in $V$. Thus, \eqref{williamsII}
also holds for every optimal alignment $v$.
If $B^n_0$ and $B^n_1$ both hold, there are at most
$(\epsilon n/16)+(\epsilon n/8)=3\epsilon n/16$ non-matched ones in each string
$X$ and $Y$ for any optimal alignment $v$. This also implies
that the number of blocks made of five zeros in $X$ which do
not satisfy the criteria to be counted towards
$N_5(v)$ is at most $2((\epsilon n/16)+(\epsilon n/8))$,
(for every optimal alignment $v$).  This implies that for every 
optimal alignment $v$, when $B^n_0$ and $B^n_1$ both hold,
then
$$N_5(v)\geq n_5-\frac{3n\epsilon}8.$$
This last inequality and \eqref{williamsII}
imply that
\begin{equation}
\label{shit}\frac{N_{5<}(v)}{n_5}\geq\left(
\frac{31}{32}- \frac{\epsilon_1}4
\right)
\left(1-\frac{3n\epsilon}{8n_5}\right),
\end{equation}
for every optimal alignment $v$.
When $B^n_3$ holds,
$$n_5\geq \left(\frac 1{32}-\frac\epsilon{16}\right)n,$$
and therefore
$$1-\frac{3n\epsilon}{8n_5}\geq
1-\frac{3\epsilon}{1/4-\epsilon/2}.$$
 Using the above inequality in \eqref{shit}
gives:
\begin{equation}
\label{shitII}\frac{N_{5<}(v)}{n_5}\geq\left(
\frac{31}{32}- \frac{\epsilon_1}4
\right)
\left(1-\frac{3\epsilon}{1/4-\epsilon/2}\right).
\end{equation}
Note that
$$\lim_{\epsilon\rightarrow 0}
\left(1-\frac{3\epsilon}{1/4-\epsilon/2}\right)=1.$$
Hence for any $\epsilon_1>0$ fixed, choosing
$\epsilon>0$ small enough, (depending on 
$\epsilon_1$), leads to
\begin{equation}
\label{conditionepsi2}
\left(
\frac{31}{32}- \frac{\epsilon_1}4
\right)
\left(1-\frac{3\epsilon}{1/4-\epsilon/2}\right)
\geq\frac{31}{32}-\frac{\epsilon_1}{2}.
\end{equation}
Assume henceforth that the inequality \eqref{conditionepsi2}
holds.  This together with \eqref{shitII} yields:
$$\frac{N_{5<}(v)}{n_5}\geq \frac{31}{32}-\frac{\epsilon_1}{2}.$$
Hence $E^n$ holds, and
$\epsilon_1>0$ and $\epsilon>0$ are chosen so that
\eqref{conditionepsi2} holds, then
$B^n_0$, $B^n_1$, $B^n_3$ and $C^n$ together imply
the event $E^n$.
\end{proof}

\medskip
Let $N_1(v)$ denote in the alignment $v$, the total number of
subsequent pairs of ones in $X$ satisfying all of the following conditions:
\begin{itemize}
\item The ones are both matched, that is aligned by $v$ with a one from $Y$.
\item Between the pair of ones in the $X$ text there is exactly one zero
and nothing else.
\item The pair of ones in the $Y$ text with which the pair from $X$ is
aligned, should contain only zeros in between them.\end{itemize}
 More precisely: for
 $v\in V^k$, let
$$N_1(v):=|\left\{
i<k: v_i=(0,0),\pi_v(i+1)-\pi_v(i)=2\right\}|,$$
where $v:=(v_1,v_2,\ldots,v_k)$.

Let $N_{1>}(v)$ denote in the alignment $v$, the total number of
subsequent pairs of ones in $X$ satisfying all of the following conditions:
\begin{itemize}
\item The ones are both matched, that is aligned by $v$ with a one from $Y$.
\item Between the pair of ones in the $X$ text there are only zeros
and, in fact, exactly one zero.
\item The pair of ones in the $Y$ text with which the pair from $X$ is
aligned, should contain only zeros in between them and,
 in fact,  two or more of them.\end{itemize}
  More precisely:
$$N_{1>}(v):=|\left\{
i<k:v_i=(0,0),\pi(i+1)-\pi(i)=2,\nu_v(i+1)-\nu_v(i)>2\right\}|.$$
Let $B^n_4$ be the event that in the sequence $X$ there are
at least $(1/4-\epsilon/16)n$ blocks of zeros of length
two.
Let $p_1(v)$ be the conditional probability on $X$, that
when picking a block of one zeros  at random in $X$,
 this block
happens to satisfy   the following two conditions:
\begin{itemize}
\item[1)] The block is contained between two consecutive matched ones.
(Matched by the alignment $v$.)
Specifically, between the two matched ones there 
is the block of length one and nothing else.
\item[2)] The pair of ones in the $Y$-text to which the pair of consecutive ones
are aligned by $v$ contains only zeros in between them and at least
two of them.\end{itemize}
In other words, $p_1(v)$ is the
conditional probability (conditional on $X$) that when picking
at random a block of zeros of length one
in $X$,  there exists
$i\leq k$, such that the randomly selected block
is equal to $[\pi_v(i-1)+1,\pi_v(i)-1]$ and all of the following conditions hold:
$$\pi_v(i)-\pi_v(i-1)=2,$$
$$\eta_v(i)-\eta_v(i-1)>2,$$
and $|v_i|=0$,
where $$v=(v_1,v_2,\ldots,v_k)\in(\mathbb{N}\times\mathbb{N})^k.$$
 Recall that the total number of blocks made of one zeros in $X$ is
denoted by $n_1$, and so the conditional probability $p_1(v)$ is equal to: 
$$p_1(v):=\frac{N_{1>}(v)}{n_1}.$$

Let us go back to the previous numerical example.
Again, let $X=10 10 10 10000011$ and let
$Y=1010001100011$.  Let $v$ be the alignment
$$
\begin{array}{c|c|c|c|c|c|c|c|c|c|c|c|c|c|c|c|c|c}
1&0&1&0& & &1&0&1&0&0&0&0&0&1 &\\\hline
1&0&1&0&0&0&1&  &1&0&0&0& & &1&1
\end{array}\,.
$$
There are three blocks of zeros of length one in $X$.
Hence, $n_1=3$.
Among such blocks, only two
are directly comprised between matched ones.  Hence,
$N_1(v)=2$.
One of the two ``suitable" blocks made of a single zero in $X$
is matched with a block of zeros in $Y$ strictly larger
than $1$.  Hence, $N_{1>}=1$ and
$p_1(v)=1/3$.

Let $F^n$ be the event
that
$$p_1(v)\geq\frac{1}{4}-\frac{\epsilon_1}{2}.
$$ for every optimal alignment $v\in V$.
Let $D^n$ be the event that
for all $v\in V$ which is an optimal alignment, we have that
$$\frac{N_{1>}(v)}{N_1(v)}\geq\frac{1}{4}- \frac{\epsilon_1}4.$$  
\begin{lemma}
Let $\epsilon>0$ and $\epsilon_1$ satisfy the
inequality \eqref{conditionepsi3} below, then
\begin{equation}
\label{B1Fn}
B^n_0\cap B^n_1\cap B^n_4\cap D^n\subset F^n,
\end{equation}
for all $n$.
\end{lemma}
\begin{proof}
 For any alignment $v\in V$, 
\begin{equation}
\label{williamsF}
\frac{N_{1>}(v)}{n_1}=\frac{N_{1>}(v)}{N_1(v)}\cdot
\frac{N_1(v)}{n_1}.
\end{equation}
Further, when the event $D^n$ holds, then, for every 
$v\in V$, 
$$\frac{N_{1>}(v)}{N_1(v)}\geq\frac{1}{4}- \frac{\epsilon_1}4.$$
Combining this last inequality with \eqref{williamsF}
yields
\begin{equation}
\label{williamsIIF}
\frac{N_{1>}(v)}{n_1}\geq\left(
\frac{1}{4}- \frac{\epsilon_1}4
\right)
\frac{N_1(v)}{n_1}.
\end{equation}
Since $B^n_0$ and $B^n_1$ together imply $B^n_2$,
then when they both hold, every
optimal alignment is in $V$. Hence, \eqref{williamsIIF}
also holds for every optimal alignment $v$.
If $B^n_0$ and $B^n_1$ both hold, there are at most
$(\epsilon n/16)+(\epsilon n/8)$ non-matched ones in each text
$X$ and $Y$ for any optimal alignment $v$. This also implies
that the number of blocks made of one zero in $X$ which do
not satisfy the criteria to be counted towards
$N_1(v)$ is at most $2[(\epsilon n/16)+(\epsilon n/8)]$
(for every optimal alignment $v$).  This implies that for every 
optimal alignment $v$, when $B^n_0$ and $B^n_1$ both hold,
then
$$N_1(v)\geq n_1-3n\epsilon/8.$$
The above inequality with \eqref{williamsIIF}
implies 
\begin{equation}
\label{shitF}\frac{N_{1<}(v)}{n_1}\geq\left(
\frac{1}{4}- \frac{\epsilon_1}4
\right)
\left(1-\frac{3n\epsilon}{8n_1}\right),
\end{equation}
for every optimal alignment $v$.
When $B^n_4$ holds,
$$n_1\geq \left(\frac 14-\frac\epsilon{16}\right)n,$$
and therefore
$$1-\frac{3n\epsilon}{8n_1}\geq
1-\frac{3\epsilon}{2-\epsilon/2}.$$
Hence, \eqref{shitF} becomes:
\begin{equation}
\label{shitIIF}\frac{N_{1>}(v)}{n_1}\geq\left(
\frac{1}{4}- \frac{\epsilon_1}4
\right)
\left(1-\frac{3\epsilon}{2-\epsilon/2}\right).
\end{equation}
Note that
$$\lim_{\epsilon\rightarrow 0}
\left(1-\frac{3\epsilon}{2-\epsilon/2}\right)=1.$$
Hence for any $\epsilon_1>0$ fixed, choosing
$\epsilon>0$ small enough, (depending on 
$\epsilon_1$), gives
\begin{equation}
\label{conditionepsi3}
\left(
\frac{1}{4}- \frac{\epsilon_1}4
\right)
\left(1-\frac{3\epsilon}{1/4-\epsilon/2}\right)
\geq\frac{31}{32}-\frac{\epsilon_1}{2}.
\end{equation}
Henceforth assume that the inequality \eqref{conditionepsi3}
holds.  This together with \eqref{shitIIF} yields:
$$\frac{N_{1>}(v)}{n_1}\geq \frac{1}{4}-\frac{\epsilon_1}{2},$$
and so $F^n$ holds. We have just proved that if
$\epsilon_1>0$ and $\epsilon>0$ are chosen so that
\eqref{conditionepsi3} holds, then
$B^n_0$, $B^n_1$, $B^n_4$ and $D^n$ together imply
the event $F^n$.
\end{proof}
\medskip

Let us return to  a numerical example.
Let $X=10101011000001$ and let
$Y=10100011000111$.  Let $v$ be the alignment
$$
\begin{array}{c|c|c|c|c|c|c|c|c|c|c|c|c|c|c|c|c|c|c}
1&0&1&0& & &1&0&1&1&0&0&0&0&0&1 & &\\\hline
1&0&1&0&0&0&1& & &1&0&0&0& & &1&1&1
\end{array}\,.
$$
Note that there are consecutive aligned ones with
one zero in between in the $X$-part.  Hence $N_1(v)=2$.
Note, that the third and fifth one in the string
$X$ are consecutive aligned ones with one zero between them.
 But there is also a non-aligned one between them,
so they do not count towards $N_1(v)$. Instead the first
consecutive aligned ones counting towards $N_1(v)$ are
given by the first and the second one in $X$. Then,
the second and third one in $X$ constitutes such
consecutive couple of ones counting towards $N_1(v)$.

Among the two pairs of consecutive aligned ones with one zero in between,
there is one which has strictly more than one zero in between the ones
in the $Y$ part.  Hence, $N_{1>}=1$.
There is one consecutive pair of aligned ones
with five zeros in between. It is given by the fifth and sixth
one in $X$.  Hence, $N_5(v)=1$.
The only pair of consecutive aligned ones with five zeros in between
in the $X$ part has 3 zeros in between in the $Y$ part. Hence,
$N_{5<}=1$.
\medskip

Recall that $A^n$ is the event that $X$ and $Y$ are such that
the inequalities \eqref{probcondA1} and
\eqref{probcondA3} are satisfied.

\begin{lemma}
Let 
\begin{equation}
\label{conditionepsi4}0<\epsilon_1<29/16,
\end{equation}
then,
\begin{equation}
\label{EFA}
E^n\cap F^n\subset A^n,
\end{equation}
holds for every $n\in \mathbb{N}$.
\end{lemma}
\begin{proof}Next,  say that a block of zeros in $X$ is {\it 
aligned with another block
of zeros} in $Y$ by $v$, if these blocks
 are in between consecutive mutually aligned ones.

Let us give an example.  
Let $X=1011$ and $Y=1001$. Let $v$ be the alignment
$$
\begin{array}{c|c|c|c|c}
1&0& &1 &1\\\hline
1&0&0&1 &
\end{array}\,.$$
In this example $X$ has one block of zeros. This block has length one.
The text $Y$ has also one block of zeros.  This block has length
two. The block of zeros in $X$ is said to get aligned by $v$ with the
block of zeros in $Y$.
Recall that to obtain $\tilde{L}_n$ from $L_n$, 
a block of five zeros is picked uniformly at random in $X$ and its length reduced by one.
Then, a zero is added to a randomly chosen block of one zero
in $X$. The modified text is denoted by $\tilde{X}$.
The optimal score between $\tilde{X}$ and $Y$ is
$\tilde{L}$.
For any optimal alignment $v$ and
when the selected block of five zeros in $X$ is aligned 
by $v$ with a block
of length strictly less than five, then the score is not reduced.
When, on top of that, the extra zero is added to
a block of length one which is aligned by $v$ to a block of at least
two zeros, than the score increases by one.  The block
of length five and the block of length one are chosen independently
from each other, and so
\begin{equation}
\label{p1p5}
\PP(\tilde{L}_n-L_n=1\mid X,Y)\geq p_5(v) p_1(v),
\end{equation} for any optimal alignment $v$.
When $E^n$ and $F^n$ both hold, then for any optimal alignment
$v$ 
\begin{equation}
\label{p1p5II}
p_5(v) p_1(v)\geq 
\left(\frac14-\frac{\epsilon_1}{2}   \right)
\left(\frac{31}{32}-\frac{\epsilon_1}{2}   \right)=
\left(\frac14\right)\left(\frac{31}{32}\right)-
\epsilon_1\left(\frac{35}{64}+\frac{\epsilon_1}{4}    \right).
\end{equation}
Since  $0<\epsilon_1<29/16$,
\eqref{p1p5} and \eqref{p1p5II}
 together imply that 
\begin{equation}
\label{A1a}
\PP(\tilde{L}_n-L_n=1\mid X,Y)\geq \left(\frac14\right)\left(\frac{31}{32}\right)-
\epsilon_1.
\end{equation}
Let $v$ be an optimal alignment. 
Since when the selected block of length five is aligned
with a block of length strictly smaller, then the score
cannot decrease, i.e.:
$$\tilde{L}_n-L_n\geq 0.$$
It follows that for any optimal alignment $v$, 
\begin{equation}
\label{tildeL-1}
\PP(\tilde{L}_n-L_n=-1\mid X,Y)\leq 1-p_5(v).
\end{equation}
Hence, when $E_n$ holds, 
\begin{equation}
\label{simple}
p_5(v)\geq \frac{31}{32}-\frac{\epsilon_1}{2}\,.
\end{equation}
Together, the inequalities \eqref{tildeL-1} and \eqref{simple}
imply
\begin{equation}
\label{A1b}
\PP(\tilde{L}_n-L_n=-1\mid X,Y)\leq \frac{1}{32}+\epsilon_1.
\end{equation}
Therefore, when $E^n$ and $F^n$ both hold,
then \eqref{A1a} and \eqref{A1b} both hold.  In other words,
$E^n$ and $F^n$ jointly imply $A^n$ when 
$\epsilon_1>0$ satisfies \eqref{conditionepsi4}.~\end{proof}

\section{Probability bounds}

\subsection{The bounds}
To begin with, let us present some useful bounds:
\begin{lemma}
There exists $\gamma_0>0$, independent of $n$, such that
\begin{equation}
\label{probbound1}\PP(B^n_0)\geq 1-e^{-\gamma_0n}.
\end{equation}
(Note that $\gamma_0$ depends  on $\epsilon$.)
\end{lemma}
\begin{proof}
The sequence $X$, resp.\ $Y$, is iid. The probability that $X_i=1$,
resp.\ $Y_i=1$ is
equal to $1/2$. Hence, by exponential inequalities, the probability
that the average $\sum_i^n X_i/n$ is different from its mean
by more than $\epsilon/16$ is exponentially small in $n$.~\end{proof}

\begin{lemma}
There exists $\gamma_3>0$, independent of $n$, such that
\begin{equation}
\label{probbound2}
\PP(B^n_3)\geq 1-e^{-\gamma_3n}.
\end{equation}
(Note that $\gamma_3$ depends on $\epsilon$.)
\end{lemma}
\begin{proof}The blocks in $X$ are iid. The probability that
a block has length five is equal to $1/32$. Again, an exponential inequality applied
to the probability that the proportion of blocks 
of zeros which has length five is below the expectation by
$\epsilon/16$. Hence that probability is exponentially small
in $n$.
\end{proof}
\medskip

Similarly,
\begin{lemma}
There exists $\gamma_4>0$, independent of $n$, such that
\begin{equation}
\label{probbound3}\PP(B^n_4)\geq 1-e^{-\gamma_4n}.
\end{equation}
(Note that $\gamma_4$ depends on $\epsilon$.)
\end{lemma}
\begin{proof}Essentially the same as the proof of the previous
lemma above.
\end{proof}
\medskip

The next lemma gives an upper bound on the number of elements
in the set $V$.
\begin{lemma}
We have that
\begin{equation}
\label{Vk}
|V^k|\leq e^{H(\epsilon/4)k}2^{\epsilon k/2},
\end{equation}
where $H$ is the entropy function.
\end{lemma}
\begin{proof}First determine which entries
are non-zero.  There are at most $\epsilon k/2$ non-zero
entries, which have to be chosen from $2k$ entries.
Hence this gives a total number of
\begin{equation}
\label{H}
\binom {2k}{\epsilon k/2}\leq e^{H(\epsilon/4)}
\end{equation}
possibilities.
Next, choose how large each entry is.
To do so, distribute among
the non-zero entries (which are already determined
in the previous step)
a total of $\epsilon k/2$ integer points. This is the same as 
finding an integer partition of the interval $2^{\epsilon k/2}$.
There are at most $2^{\epsilon k/2}$ integer partitions
of the interval $[0,\epsilon k/2]$.  This, when combined, with \eqref{H} indicates
that $V^k$ contains no more than
$e^{H(\epsilon/4)}2^{\epsilon k/2}$
elements.~\end{proof}

\medskip

Eventually, we have
\begin{lemma}
Let $\epsilon>0$ and $\epsilon_1$ be such that
\eqref{conditionepsi5} and \eqref{conditionepsi6} both hold,
then 
\begin{equation}
\label{probbound4}
\PP(C^n)\geq 1-e^{-\gamma_cn},
\end{equation}
where $\gamma_c$ is a positive constant independent of $n$,
(but depending on $\epsilon$).
\end{lemma}

\begin{proof} Let $C^n_*$ be the event that
if $N_5(v)\geq n/33$, then $C^n$ holds.
In other words
$$C^n_*=\left(C^n\cap 
\left\{N_5(v)\geq \frac{n}{33}\right\}\right)\cup
\left\{N_5(v)< \frac{n}{33}\right\}.$$
When $B^n_3$ holds there are at least 
$(1/32-\epsilon/16)n$ blocks of five zeros in $X$.
When $B^n_0$ and $B^n_1$ both hold, then
it was argued that there are at most $6\epsilon n/16$ ones
not matched total in both sequences $X$ and $Y$ for any optimal alignment
$v$.

Hence, when $B^n_0$, $B^n_1$ and $B^n_3$ all hold,
then for any optimal alignment $v$, 
\begin{equation}
\label{little}
N_5(v)\geq \left(\frac1{32}-\frac{7\epsilon}{16}\right)n.
\end{equation}
For $\epsilon>0$ small enough, 
\begin{equation}
\label{conditionepsi5}
\frac{1}{32}-\frac{7\epsilon }{16}\geq \frac{1}{33}.
\end{equation}
From here on, assume that \eqref{conditionepsi5} holds,
so that \eqref{little} implies
that
$$N_5(v)\geq\frac{n}{33}.$$
Hence
$$B^n_0\cap B^n_1\cap B^n_3\subset \left\{N_5(v)\geq\frac{n}{33}
\right\},$$
and therefore
\begin{equation}
\label{useful}
\left\{N_5(v)<\frac{n}{33}
\right\}\subset (B^{n}_0)^c\cup (B^{n}_1)^c \cup (B^{n}_3)^c.
\end{equation}  
But
$$(C^{n})^c\subset (C^{n}_*)^c\cup \left\{N_5(v)<\frac{n}{33}
\right\},$$
and so with the help of \eqref{useful} and of Lemma~\ref{lem4.1}, 
$$(C^{n})^c\subset (C^{n}_*)^c\cup (B^{n}_0)^c \cup (B^{n}_3)^c.$$
The last inclusion implies that
$$\PP((C^{n})^c)\leq \PP((C^{n}_*)^c)+\PP((B^{n}_0 )^c)+ 
\PP((B^{n}_3)^c).$$
We already proved exponentially small upper bounds for
$\PP((B^{n}_0)^c)$ and for $\PP((B^{n}_3)^c)$.  Hence it only remains
to prove a similar upper bound for 
$\PP((C^{n}_*)^c)$.
Let $Z_1,Z_2,\ldots$ be a sequence of iid geometric random 
variables with parameter $1/2$. Let $W_i$ be the indicator variable
which is equal to $1$ if $Z_i<5$. Then, $\PP(W_i=1)=31/32$, and
\begin{equation}
\label{W1}
\PP((C^{n}_*)^c)\leq \sum_{k=n/33}^{\infty}|V^k|
\PP\left(\frac{W_1+W_2+\cdots+W_{k}}{k}<\frac{31}{33}-
\frac{\epsilon_1}4\right).
\end{equation}
Classically
\begin{equation}
\label{largedevW}
\PP\left(\frac{W_1+W_2+\cdots+W_{n/33}}{n/33}<\frac{31}{33}-
\frac{\epsilon_1}4\right)
\leq e^{-n\gamma(\epsilon_1)},
\end{equation}
for some constant $\gamma(\epsilon_1)$.
Combining \eqref{Vk} and \eqref{largedevW}, leads to
\begin{equation}
\label{maindeviation}
|V^k|
\PP\left(\frac{W_1+W_2+\ldots+W_{k}}{k}<\frac{31}{33}-
\frac{\epsilon_1}4\right)
\leq e^{H(\epsilon)k}2^{\epsilon k} e^{-k\gamma(\epsilon_1)}.
\end{equation}
Hence \eqref{W1} and \eqref{maindeviation} together provide
an exponential upper bound for $\PP((C^{n}_*)^c)$ as soon as the
following  inequality:
\begin{equation}\label{conditionepsi6}
H(\epsilon)+\epsilon-\gamma(\epsilon_1)<0
\end{equation}
is satisfied.
Note that
$$\lim_{\epsilon\rightarrow 0}\left(H(\epsilon)+\epsilon\right)=0,$$
while $\gamma(\epsilon_1)>0$, for every $\epsilon_1>0$.
This implies that for any $\epsilon_1>0$ fixed, 
\eqref{conditionepsi6} holds, when taking
$\epsilon>0$ small enough.
\end{proof}
\medskip

In complete similarly to the previous lemma we have:

\begin{lemma}
Let $\epsilon>0$ and $\epsilon_1$ be such that
\eqref{conditionepsi5} and \eqref{conditionepsi6} both hold,
then there exists $\gamma_d>0$ independent of
$n$ such that
\begin{equation}
\label{probbound5}
\PP(D^n)\geq 1-e^{-\gamma_dn}.
\end{equation}
\end{lemma}
\begin{proof}This proof is similar to the previous one and so is omitted. 
\end{proof}

\subsection{Overview}\label{overview}
In Section~3 we proved that ${\rm Var\;}{L_n}$ is of order $n$.
For this we assumed that Theorem~2.2 holds.
So it still remains to prove Theorem~2.2.
First, let us mention that the order $n$ for
${\rm Var\; }L_n$ follows from two things:
\begin{itemize}
\item[{\bf a)}] $\Delta_n:=\tilde{L}_n-L_n$ needs to have a positive bias.
More precisely, we want that for any $(x,y)\in A^n$ and any $n$,
$$\E\left(\tilde{L}_n-L_n\mid X=x,Y=y\right)>K,$$
where $K>0$ is some positive constant.
\item[{\bf b)}] The probability $\PP(A^n)$ needs to be close to one,
more precisely $1-\PP(A^n)$ needs to be no more than a
stretched negative exponential.\end{itemize}
Let us first mention problem a). Note that between
$\tilde{L}_n$ and $L_n$ the score can change by at most one.
In other words,
\begin{equation}
\label{one}
\PP(\tilde{L}_n-L_n\in\{-1,0,1\})=1.
\end{equation}
From the equation \eqref{one} and with the help of \eqref{probcondA1}
and \eqref{probcondA3}, 
\begin{equation}
\label{expectationDL}
\E\left(\tilde{L}_n-L_n= 1\mid X=x,Y=y\right)\geq \frac{31}{128}-\frac{1}{32}
- 2\epsilon_1=\frac{27}{128}-2\epsilon_1.
\end{equation}
Taking $\epsilon_1>0$ small enough so that
\begin{equation}
\label{conditionbias}
\epsilon_1<\frac{3^3}{2^8},
\end{equation}
ensures the positive bias.

Let us next discuss problem b). 
In the previous section, the following inclusions were shown to
hold:
\begin{itemize}
\item $B_0^n\subset B_1^n$, when \eqref{conditionepsi0} holds.
\item $B_0^n\cap B_1^n\subset B_2^n$, when \eqref{conditionepsi1} holds.
\item $B_0^n\cap B_1^n\cap B^n_3\cap C^n\subset E^n$, when 
\eqref{conditionepsi2} holds.
\item $B_0^n\cap B_1^n\cap B^n_4\cap D^n\subset F^n$, when 
\eqref{conditionepsi3} holds.
\item $E^n\cap F^n\subset A^n$, 
when \eqref{conditionepsi4} holds.\end{itemize}
These inclusions imply, when all the conditions
\eqref{conditionepsi0}, \eqref{conditionepsi1}, \eqref{conditionepsi2},
\eqref{conditionepsi3} and
\eqref{conditionepsi4} hold, 
$$ B^n_0\cap B^n_3\cap B^n_4\cap E^n\cap D^n\subset A^n,$$
and therefore
\begin{equation}
\label{mainoin}
\PP((A^{n})^c)\leq 
\PP((B^{n}_0)^c)+\PP((B^{n}_3)^c)+\PP((B^{n}_4)^c)+\PP((C^{n})^c)
+\PP((D^{n})^c).
\end{equation}
The inequality \eqref{mainoin} implies that $\PP((A^{n})^c)$ is 
exponentially small in $n$, as soon
as exponential bounds are available
for $\PP((B^{n}_0)^c)$, $\PP((B^{n}_3)^c)$, $\PP((B^{n}_4)^c)$, 
$\PP((C^{n})^c)$ and $\PP((D^{n})^c)$.
The probabilities $\PP((B^{n}_0)^c)$, $\PP((B^{n}_3)^c)$ and 
$\PP((B^{n}_4)^c)$,
only depend on $\epsilon$.  For any $\epsilon>0$,
the inequalities \eqref{probbound1}, \eqref{probbound2} and
\eqref{probbound3} provide exponentially small bounds
for $\PP((B^{n}_0)^c)$, $\PP((B^{n}_3)^c)$ and $\PP((B^{n}_4)^c)$.  Hence,
no special condition on $\epsilon>0$ and
$\epsilon_1>0$ are needed in order to ensure that 
$\PP((B^{n}_0)^c)$, $\PP((B^{n}_3)^c)$ and $\PP((B^{n}_4)^c)$
are exponentially small in $n$.

The inequalities \eqref{probbound4} and \eqref{probbound5}
also provide exponentially small upper bounds are also obtained
for $\PP((C^{n})^c)$ and $\PP((D^{n})^c)$. 
However these bounds only hold if $\epsilon$ and $\epsilon_1$ satisfy
\eqref{conditionepsi5} and
\eqref{conditionepsi6}.

To prove that $\PP((A^{n})^c)$ is exponentially small it thus
remains to prove that there exists $\epsilon$ and $\epsilon_1$
satisfying all the following:\begin{itemize}
\item[a)] The condition \eqref{conditionbias}, which ensures the 
conditional bias \eqref{expectationDL} on
$\tilde{L}_n-L_n$.
\item[b)] All the conditions for the inclusions.  These are the inequalities 
\eqref{conditionepsi1},
\eqref{conditionepsi2} and
\eqref{conditionepsi3}.
\item[c)] The conditions for the exponentially small upper bounds 
for $\PP((C^{n})^c)$ and $\PP((D^{n})^c)$. These are the inequalities 
\eqref{conditionepsi5} and
\eqref{conditionepsi6}.
\end{itemize}
To see that there exist $\epsilon,\epsilon_1>0$ satisfying
all the above conditions simultaneously, note that these conditions
can be classified into three types:\begin{description}
\item[Type I:] conditions involving only $\epsilon_1$.  These conditions
all hold for $\epsilon_1>0$ small enough.
\item[Type II:] conditions involving $\epsilon_1$ and $\epsilon$.
All these conditions are such that for any $\epsilon_1>0$
fixed, they hold as soon as $\epsilon>0$ is small enough.
\item[Type III:] conditions involving only $\epsilon$. They all hold
as soon as $\epsilon$ is taken small enough.\end{description}

It is now easy to see that there exists $\epsilon_1>0$ and
$\epsilon>0$ such that all the conditions 
\eqref{conditionepsi0},
\eqref{conditionepsi1},
\eqref{conditionepsi2},
\eqref{conditionepsi3},
\eqref{conditionepsi4},
\eqref{conditionepsi5},
\eqref{conditionepsi6},
 and
\eqref{conditionbias} 
simultaneously hold. For this choose first
$\epsilon_1>0$ small enough so that all equations of type I
are satisfied.  Then choose $\epsilon>0$ small enough so that
all conditions of type II and type III are satisfied.
Summarizing:
There exist $\epsilon,\epsilon_1>0$ and $c_1,s^*>0$
not depending on $n$, such that if $s(1,1)\geq s^*$, then
\begin{equation}\label{finalfinal}
\PP((A^{n})^c)\leq e^{c_1n},
\end{equation} 
for all $n$.

\bibliographystyle{plain}
\bibliography{bio15}

\end{document}